\documentclass[11pt]{amsart}
\usepackage[margin=1.15in]{geometry}
\usepackage{amsmath,amssymb,amsthm,mathtools}
\usepackage[protrusion=true,expansion=false]{microtype}
\usepackage{enumitem,booktabs}
\usepackage[colorlinks=true,linkcolor=blue!50!black,citecolor=blue!50!black,urlcolor=blue!50!black]{hyperref}
\usepackage{xcolor}

\setlist{itemsep=3pt,topsep=5pt}
\allowdisplaybreaks

\newtheoremstyle{airy}{11pt}{11pt}{\itshape}{}{\bfseries}{.}{0.6em}{}
\newtheoremstyle{airydef}{11pt}{11pt}{\normalfont}{}{\bfseries}{.}{0.6em}{}
\newtheoremstyle{airyrem}{9pt}{9pt}{\normalfont}{}{\itshape}{.}{0.6em}{}

\theoremstyle{airy}
\newtheorem{mainthm}{Theorem}

\newtheorem*{thmnn}{Theorem}
\newtheorem{theorem}{Theorem}[section]
\newtheorem{lemma}[theorem]{Lemma}
\newtheorem{proposition}[theorem]{Proposition}
\newtheorem{corollary}[theorem]{Corollary}

\theoremstyle{airydef}
\newtheorem{definition}[theorem]{Definition}

\theoremstyle{airyrem}
\newtheorem{remark}[theorem]{Remark}

\newcommand{\QQ}{\mathbb{Q}}
\newcommand{\ZZ}{\mathbb{Z}}
\newcommand{\CC}{\mathbb{C}}
\newcommand{\FF}{\mathbb{F}}

\newcommand{\Kbar}{\overline{K}}
\newcommand{\hh}{\widehat{h}}
\newcommand{\OO}{\mathcal{O}}

\newcommand{\p}{\mathfrak{p}}
\DeclareMathOperator{\Res}{Res}

\begin{document}

\title[An elliptic Lehmer bound of exponent one]{An elliptic Lehmer bound of exponent one at a supersingular prime}
\author{Jorge Mello}
\date{Working note, November 2026}

\begin{abstract}
Let $E$ be an elliptic curve over a number field $K$ with good supersingular reduction at a prime $\p$. We prove that
the N\'eron--Tate height of every non-torsion $P\in E(\Kbar)$ satisfies $\hh_E(P)\ge c\,D^{-1}\exp(-c'\sqrt{\log3D})$,
where $D=[K(x(P)):K]$ and $c,c'>0$ depend only on $E$, $K$ and $\p$; so
$\hh_E(P)\ge D^{-1-o(1)}$. By a result of Elkies this applies to every elliptic curve over $\QQ$, and it improves the exponents of
Masser and of Hindry and Silverman.
\end{abstract}

\maketitle

\section{Introduction}

Let $E$ be an elliptic curve over a number field $K$ and denote by $\hh_E$ its N\'eron--Tate height. For a non-torsion point $P$
of degree $D=[K(P):K]$, the elliptic Lehmer conjecture asks for a real constant $c(E/K)>0$ such that $$\hh_E(P)\ge \dfrac{c(E/K)}{D}.$$ For curves with complex multiplication,
Laurent \cite{Lau} proved the Dobrowolski-type bound $c\,D^{-1}(\log\log3D/\log3D)^3$, using Frobenius endomorphisms at many split
primes. For curves with non-integral $j$-invariant, Hindry and Silverman \cite{HSi} proved $c/(D^2\log^23D)$, and for all curves
Masser \cite{Mas} proved $c/(D^3\log^23D)$; for integral $j$-invariant and without complex multiplication this was, as far as we
know, the best exponent.

Throughout, $K$ is a number field and $\p$ a prime of $K$ above the rational prime $p$, with ramification index $e$
and residue degree $f$. Let $\OO=\OO_{K,\p}$ be the local ring, $k=\OO/\p$ and $N\p=p^f$. Let $E/K$ be given by a
Weierstrass equation with coefficients in $\OO$ and discriminant in $\OO^\times$ whose reduction $\widetilde E$ is
supersingular; every elliptic curve with good supersingular reduction at $\p$ has such an equation. For $P\in E(\Kbar)$ let
$e_\p(P)$ be the largest absolute ramification index of a prime of $K(P)$ above $\p$, so that $e_\p(P)\le e[K(P):K]$.

Our main result is a lower bound with exponent $1+o(1)$ for every curve with a prime of good supersingular reduction. It
improves Masser's exponent $3$ for all such curves and, when $j(E)$ is not integral, the exponent $2$ of Hindry and Silverman.

\begin{mainthm}\label{thm:main}
There are $c_1,c_2>0$, depending only on $E$, $K$ and $\p$, such that every non-torsion $P\in E(\Kbar)$ satisfies, with
$D=[K(x(P)):K]$,
\begin{enumerate}[label=\rm(\roman*),leftmargin=*]
\item $\displaystyle\hh_E(P)\ \ge\ \frac{c_1}{e_\p(P)\,D}\exp\bigl(-c_2\sqrt{\log3D}\bigr)$;
\item $\displaystyle\hh_E(P)\ \ge\ \frac{c_1}{D}\exp\bigl(-(c_2+1)\sqrt{\log3D}\bigr)$.
\end{enumerate}
In particular $\hh_E(P)\ge D^{-1-o(1)}$ for all non-torsion $P$.
\end{mainthm}

By a result of Elkies, the hypothesis holds for every elliptic curve over $\QQ$, and more generally in the following cases.

\begin{corollary}\label{cor:Q}
Let $E$ be an elliptic curve over $\QQ$, or more generally over a number field $K$ such that $E$ has complex
multiplication or $\QQ(j(E))$ has a real embedding. Then there are $c,c'>0$, depending only on $E$ and $K$, such that every non-torsion
$P\in E(\Kbar)$ satisfies, with $D=[K(x(P)):K]$,
$$\hh_E(P)\ \ge\ \frac cD\exp\bigl(-c'\sqrt{\log3D}\bigr);$$
in particular $\hh_E(P)\ge D^{-1-o(1)}$.
\end{corollary}

\begin{proof}
Such a curve has a prime of good supersingular reduction: over $\QQ$ by Elkies \cite{Elk87}, and in general by
Proposition \ref{prop:ss}. Apply Theorem \ref{thm:main}.
\end{proof}

Since $[K(x(P)):K]\le[K(P):K]$ and the bounds decrease with $D$, they also hold with $D=[K(P):K]$, and the exponent $1$ is
optimal: if $Q$ is non-torsion over a finite extension $K'$ of $K$ and $[m]P=Q$, then $[K(P):K]\le m^2[K':K]$ and
$\hh_E(P)=\hh_E(Q)/m^2$ (see also Remark \ref{rem:resdegsharp}). For curves with complex multiplication Laurent's bound is stronger. Part (i) is the better bound when $e_\p(P)<\exp(\sqrt{\log3D})$.

Here, among other things, we adapt Dobrowolski's method, which for $x\mapsto x^p$ uses the congruence
$A(x^p)\equiv A(x)^p\pmod p$ for many primes $p$, to the single prime $\p$. Let $\varphi$ be the Latt\`es map of $[p]$, the rational map of degree $p^2$ with
$x([p]Q)=\varphi(x(Q))$ for all $Q\in E$. Its reduction modulo $\p$ is $x\mapsto(x^{p^2}+s)/c$ for some $c\ne0$ and $s$ (Lemma
\ref{lem:redp}), and we deduce from this that the binary forms $A_n$ whose zeros are the conjugates of $\alpha_n=x([p^n]P)$
satisfy $A_{n+q}\equiv\pm A_n\pmod{\p^{n+1}}$ for a fixed $q\ge1$ (Theorem \ref{thm:cong}). This is combined with Baker's lower
bound for the energy of the dynamical Green functions of $\varphi$, that is, for sums of the Green function of $\varphi$ over pairs of
points, along $\asymp\sqrt{\log D}$ consecutive steps of the orbit. Write $D_n=[K(\alpha_n):K]$ and
$r_n=D_n/D_{n+1}\in\{1,\dots,p^2\}$; the degree \emph{drops} at step $n$ if $r_n\ge2$. We obtain the following.

\begin{thmnn}
There are $c,c'>0$ and $k(D)=O(\sqrt{\log3D})$ such that every non-torsion $P\in E(\Kbar)$ satisfies the following: if
$D_n$ is constant for $m\le n\le m+k(D_m)$, then
$$\hh_E(P)\ \ge\ \frac cD\exp\bigl(-c'\sqrt{\log3D}\bigr)\prod_{n<m}\frac{r_n}{p^2}.$$
\end{thmnn}

This is Theorem \ref{thm:window} applied to the Latt\`es map, by Lemma \ref{lem:hE}. With $m=0$ it gives
$\hh_E(P)\ge c\,D^{-1}\exp(-c'\sqrt{\log3D})$ whenever the degree does not drop during the first $k(D)$ steps, for instance
when every prime factor of $D$ exceeds $p^2$, since each $r_n$ is an integer in $[1,p^2]$ dividing $D$. In general every step
$n<m$ with $r_n<p^2$ lowers the bound by the factor $r_n/p^2$, and one needs to control where the degree drops.

To this end we study the dynamics of $\varphi$ at $\p$. A point $z$ is \emph{periodic} if $\varphi^n(z)=z$ for some $n\ge1$, and
\emph{critical} if $\varphi$ is not injective near $z$; the Latt\`es map has no periodic critical point (Lemma \ref{lem:crit}). Each
residue disc of $\mathbb P^1(\CC_p)$, the set of points with a given reduction modulo $\p$, contains exactly one periodic point.
If the degree drops at step $n$, two conjugates of $\alpha_n$ have the same image under $\varphi$. If both lie above conjugates of
$x(P)$ whose fields over $K_\p$ have small ramification index, they are $\p$-adically close to the periodic point $\eta$ of their
residue disc, and a divided difference shows that the derivative of $\varphi$ at $\eta$ has large $\p$-adic valuation. This is
impossible, because this valuation is bounded on the periodic points, none of which is critical. Applying this to every fibre of $\varphi$, we obtain that a drop at step $m$ forces at least half of the conjugates of
$x(P)$ to have ramification index $\gg p^{2m}$ above $\p$ (Lemma \ref{lem:half}). Taking $m$ to be the last drop gives part (i).
For part (ii) we use that these ramified conjugates form $\p$-adic clusters, sets of $\asymp e_\p(P)$ conjugates that are pairwise
$\p$-adically close, so that the conjugates of $\alpha_m$ have large energy at $\p$, that is, a large sum of $-\log$ of their
pairwise $\p$-adic chordal distances (Proposition \ref{prop:energy}). This gives $\hh_E(P)\gg D^{-1}$ unless $D_mp^{2m}\ll D\log D$, and in that case the bound
above, applied after the last drop, gives $\hh_E(P)\gg(D_mp^{2m})^{-1}\exp(-c'\sqrt{\log3D})$.

As a byproduct of the local analysis, we obtain a bound in terms of the ramification alone when the residue degrees above
$\p$ are bounded: the factor $\exp(-c_2\sqrt{\log3D})$ of Theorem \ref{thm:main} disappears, and the dependence on
the ramification is sharp.

\begin{mainthm}\label{thm:resdeg}
Let $f_0\ge1$. There is $c_3>0$, depending only on $E$, $K$, $\p$ and $f_0$, such that every non-torsion $P\in E(\Kbar)$ for which
every prime of $K(P)$ above $\p$ has absolute residue degree at most $f_0$ satisfies
$$\hh_E(P)\ \ge\ \frac{c_3}{e_\p(P)}\ \ge\ \frac{c_3}{2e\,D},\qquad D=[K(x(P)):K].$$
In particular $\hh_E(P)$ is bounded below by a positive constant if moreover the ramification above $\p$ is bounded.
\end{mainthm}

The exponent of $e_\p(P)$ cannot be lowered (Remark \ref{rem:resdegsharp}). The proof (\S\ref{sec:resdeg}) compares a single pair of
levels of the orbit; it uses that $[p]$ contracts the formal group $\p$-adically near each residue class (a fibre of the reduction map
$E(\CC_p)\to\widetilde E(\overline{\FF}_p)$), and that the conjugates of
$P$ lie in a bounded number of residue classes. For the Weil height, Bombieri and Zannier \cite{BZ} gave lower bounds over fields of
bounded local degree in which the ramification index enters in the same way, and Baker and Petsche \cite{BP} proved the
Bogomolov property for elliptic curves over totally $p$-adic fields.

We summarise in Table \ref{tab:history} the bounds above together with ours, in the form $\hh_E(P)\ge D^{-\ell-o(1)}$, where
$D=[K(P):K]$ and $o(1)$ denotes a quantity that tends to $0$ as $D\to\infty$ and depends only on $E$, $K$ (and $\p$, $f_0$).

\begin{table}[h!]
\centering\small
\renewcommand{\arraystretch}{1.25}
\begin{tabular}{|c|l|l|}
\hline
\rule[-6pt]{0pt}{18pt}Lower bound for $\hh_E(P)$ & Restriction on $E/K$ or $P$ & Reference\\
\hline
\rule[-6pt]{0pt}{18pt}$D^{-1-o(1)}$ & $E$ has complex multiplication & Laurent 1983 \cite{Lau}\\
\hline
\rule[-6pt]{0pt}{18pt}$D^{-3-o(1)}$ & none & Masser 1989 \cite{Mas}\\
\hline
\rule[-6pt]{0pt}{18pt}$D^{-2-o(1)}$ & $j(E)$ non-integral & Hindry--Silverman 1990 \cite{HSi}\\
\hline
$D^{-1-o(1)}$ & {\renewcommand{\arraystretch}{1.1}\begin{tabular}{@{}l@{}}\rule{0pt}{11pt}$E$ has a prime of good supersingular\\ reduction (e.g.\ every $E/\QQ$)\rule[-5pt]{0pt}{0pt}\end{tabular}} & Theorem \ref{thm:main}\\
\hline
$c\,e_\p(P)^{-1}\ge c(2eD)^{-1}$ & {\renewcommand{\arraystretch}{1.1}\begin{tabular}{@{}l@{}}\rule{0pt}{11pt}$E$ supersingular at $\p$, residue degrees\\ of $K(P)$ above $\p$ at most $f_0$\rule[-5pt]{0pt}{0pt}\end{tabular}} & Theorem \ref{thm:resdeg}\\
\hline
\end{tabular}
\vspace{4pt}
\caption{Lower bounds in the elliptic Lehmer problem.}\label{tab:history}
\end{table}

Some of the results proven here can also be used to obtain improvements for abelian varieties.
Let $A$ be an abelian variety of dimension $g$ over $K$, $\mathcal L$ an ample symmetric line bundle on $A$, and
$\hh_{A,\mathcal L}$ the N\'eron--Tate height. For non-torsion $P\in A(\Kbar)$ of degree $D=[K(P):K]$, Masser proved
$\hh_{A,\mathcal L}(P)\ge c\,D^{-(2g+1)}(\log3D)^{-2g}$ \cite{Mas2,GR}. Lastly, we deduce from Theorem \ref{thm:main} the
exponent $1+o(1)$ for abelian varieties isogenous to products of elliptic curves, in place of Masser's $2g+1$.

\begin{mainthm}\label{thm:abelian}
Let $K'/K$ be a finite extension, and suppose that $A_{K'}$ is isogenous over $K'$ to a product $E_1\times\dots\times E_g$
of elliptic curves over $K'$, each of which has a prime of good supersingular reduction. Then there are $c,c'>0$, depending only on $A$, $\mathcal L$, $K'$ and the isogeny, such that every
non-torsion $P\in A(\Kbar)$ satisfies
$$\hh_{A,\mathcal L}(P)\ \ge\ \frac cD\exp\bigl(-c'\sqrt{\log3D}\bigr),\qquad D=[K(P):K].$$
\end{mainthm}

The hypothesis holds for every abelian variety over $\QQ$ isogenous over $\QQ$ to a product of elliptic curves, and for
Weil restrictions of elliptic curves with complex multiplication or with $\QQ(j(E))$ having a real embedding, some of which are simple over
$\QQ$ (Corollary \ref{cor:abex}).

In Sections \ref{sec:setting}--\ref{sec:green} we prove the bound for a rational map of degree $p^2$ that lifts, up to an affine change of variable, the
$p^2$-power Frobenius, with no use of elliptic curves. Section \ref{sec:redp} contains the needed background on the Latt\`es
map of $[p]$, and Section \ref{sec:local} its dynamics at $\p$. Theorem \ref{thm:main} is proved in \S\ref{sec:clusters},
Theorem \ref{thm:resdeg} in \S\ref{sec:resdeg} and Theorem \ref{thm:abelian} in \S\ref{sec:abelian}, and \S\ref{sec:ss} discusses
supersingular primes.

\section{Maps that lift Frobenius}\label{sec:setting}

In this section we fix a rational map of degree $p^2$ whose reduction at $\p$ is a Frobenius followed by an affine map, and state a bound that works when there is a set of consecutive iterates of constant degree. This is proved in \S\S\ref{sec:forms}--\ref{sec:green}.

Let $K$, $\p$, $p$, $e$, $f$, $\OO$ and $k$ be as in the introduction, and put $d=p^2$. In \S\S\ref{sec:setting}--\ref{sec:green} we fix a pair
$\Phi=(\Phi_0,\Phi_1)$ of binary forms of degree $d$ with coefficients in $\OO$, and elements $c\in\OO^\times$, $s\in\OO$, such that
\begin{equation}\label{eq:H}
\Phi\ \equiv\ \Psi\pmod\p,\qquad\Psi=\bigl(X^d+sZ^d,\ cZ^d\bigr),
\end{equation}
the congruence being coefficientwise. The resultant $\Res(\Phi)=\Res(\Phi_0,\Phi_1)$ is a unit in $\OO$, since its reduction is
$\Res(X^d+sZ^d,cZ^d)=\pm c^d\ne0$. So $\Phi_0,\Phi_1$ have no common zero in $\mathbb P^1(\Kbar)$, and $\Phi$ defines a rational map
$\varphi\colon\mathbb P^1\to\mathbb P^1$ of degree $d$ over $K$, whose reduction modulo $\p$ is $x\mapsto(x^{p^2}+\bar s)/\bar c$.

Let $h$ be the absolute logarithmic Weil height on $\mathbb P^1(\Kbar)$ and
$$\hh_\varphi(\beta)=\lim_{n\to\infty}d^{-n}h\bigl(\varphi^n(\beta)\bigr)$$
the canonical height of $\varphi$ \cite{CS}. It satisfies $\hh_\varphi\circ\varphi=d\,\hh_\varphi$, it is invariant under
$\mathrm{Gal}(\Kbar/K)$ because $\varphi$ is defined over $K$, and $\hh_\varphi(\beta)=0$ exactly when $\beta$ is preperiodic, by Northcott's
theorem. For $\alpha\in\mathbb P^1(\Kbar)$ we write
$$\alpha_n=\varphi^n(\alpha),\qquad D_n=[K(\alpha_n):K],\qquad r_n=D_n/D_{n+1}.$$
Since $K(\alpha_{n+1})\subseteq K(\alpha_n)$, each $r_n$ is an integer, and $r_n\le d$ because $\alpha_n$ lies in $\varphi^{-1}(\alpha_{n+1})$, a
set of at most $d$ points stable under $\mathrm{Gal}(\Kbar/K(\alpha_{n+1}))$.

The following bound, proved in \S\ref{sec:green}, bounds $\hh_\varphi(\alpha)$ from below as soon as the degree is constant along $\asymp\sqrt{\log D}$ consecutive steps of the orbit.

\begin{theorem}\label{thm:window}
There are constants $c_1',c_2>0$ and a non-decreasing function $k(D)=k_0+\lceil B\sqrt{\log3D}\,\rceil$, with $k_0$, $B$
depending only on $\Phi$, $K$ and $\p$, such that the following holds. Let $\alpha\in\mathbb P^1(\Kbar)$ be not preperiodic for
$\varphi$, and let $m\ge0$ be such that $D_{m+i}=D_m$ for $0\le i\le k(D_m)$. Then
$$\hh_\varphi(\alpha)\ \ge\ \mathcal B(D_0)\prod_{n<m}\frac{r_n}d,\qquad\mathcal B(D)=\frac{c_1'}D\exp\bigl(-c_2\sqrt{\log3D}\bigr).$$
In particular, if $D_n=D_0$ for $0\le n\le k(D_0)$, then $\hh_\varphi(\alpha)\ge\mathcal B(D_0)$; and steps with $r_n=d$ do not lower
the bound.
\end{theorem}

The constants are made explicit in \S\ref{sec:green}. The proof occupies \S\S\ref{sec:forms}--\ref{sec:green}: we construct
the characteristic forms $A_n$ of the orbit, whose zeros are the $K$-conjugates of $\alpha_n$, prove the congruence $A_{n+q}\equiv\pm A_n\pmod{\p^{n+1}}$, and compare it with an
archimedean estimate.

\section{Binary forms and push-forward}\label{sec:forms}

In this section we recall binary forms and resultants, and define the push-forward of a form by a pair of forms, whose zeros are the images of the zeros of the form.

Let $R$ be a commutative ring. A binary form of degree $D$ over $R$ is $A=\sum_{i=0}^Da_iX^iZ^{D-i}\in R[X,Z]$. For an
ideal $I$ of $R$ we write $A\equiv B\pmod I$ when all coefficients of $A-B$ lie in $I$. If $R=\OO$ is a discrete
valuation ring, $A$ is \emph{primitive} if some coefficient is a unit; by Gauss's lemma a product of primitive forms is
primitive.

For forms $A$ of degree $D$ and $G$ of degree $\delta$, $\Res(A,G)$ denotes the Sylvester resultant. It is a polynomial with
integer coefficients in the coefficients of $A$ and $G$, homogeneous of degree $\delta$ in those of $A$ and of degree $D$ in
those of $G$, and multiplicative in each argument. Over a field, if $A=\prod_{i=1}^D(z_iX-x_iZ)$, then
\begin{equation}\label{eq:resprod}
\Res(A,G)=\pm\prod_{i=1}^DG(x_i,z_i),
\end{equation}
with a sign depending only on $D$ and $\delta$. Every non-zero form over an algebraically closed field factors in this way.

The push-forward of a form $A$ by a pair $\Phi$ is a form whose zeros are the images under $\Phi$ of the zeros of $A$.

\begin{definition}\label{def:push}
Let $\Phi=(\Phi_0,\Phi_1)$ be a pair of binary forms of degree $\delta$ over $R$. For a form $A$ of degree $D$ over $R$ put
$$\Phi_*A(X,Z)\ =\ \Res_{(x,z)}\bigl(A(x,z),\ X\,\Phi_1(x,z)-Z\,\Phi_0(x,z)\bigr),$$
the resultant being taken in the variables $(x,z)$ with coefficients in $R[X,Z]$.
\end{definition}

Since the second argument is linear in $(X,Z)$, $\Phi_*A$ is a form of degree $D$. Its coefficients are polynomials with
integer coefficients in the coefficients of $A$ and $\Phi$, and $\Phi_*(\lambda A)=\lambda^\delta\Phi_*A$. Over a field,
\eqref{eq:resprod} gives
\begin{equation}\label{eq:pushprod}
\Phi_*\prod_{i=1}^D(z_iX-x_iZ)\ =\ \pm\prod_{i=1}^D\bigl(\Phi_1(x_i,z_i)X-\Phi_0(x_i,z_i)Z\bigr).
\end{equation}
So if $\Phi_0,\Phi_1$ have no common zero in $\mathbb P^1$, then $\Phi_*A$ vanishes exactly at the images under $\Phi$ of the
zeros of $A$, with multiplicities.

The next lemma turns a congruence between two forms into a divisibility of their resultant.

\begin{lemma}\label{lem:resmod}
Let $\OO$ be a discrete valuation ring and $I$ an ideal of $\OO$. Let $A,B\in\OO[X,Z]$ have degree $D$ with $A\ne0$, and
suppose $B\equiv\mu A\pmod{I^m}$ with $\mu\in\OO$. Then $\Res(A,B)\in I^{mD}$.
\end{lemma}

\begin{proof}
Write $B=\mu A+C$ with the coefficients of $C$ in $I^m$. Over an algebraic closure of the fraction field write
$A=\prod_i(z_iX-x_iZ)$. By \eqref{eq:resprod}, since $A(x_i,z_i)=0$, $\Res(A,B)=\pm\prod_iB(x_i,z_i)=\pm\prod_iC(x_i,z_i)
=\Res(A,C)$, the signs being the same because the degrees are. As $\Res(A,C)$ is homogeneous of degree $D$ in the
coefficients of $C$ with coefficients in $\OO$, it lies in $I^{mD}$.
\end{proof}

\section{Congruences for push-forwards}\label{sec:cong}

In this section we show that push-forward by a pair of forms that lifts Frobenius raises congruences between forms by one power of $\p$.

Throughout this section $R$ is a commutative ring and $I\subseteq R$ an ideal with $p\in I$; for $m\ge0$ we write
$A\equiv B\pmod{I^m}$ when all coefficients of $A-B$ lie in $I^m$ (with $I^0=R$). Let $\zeta$ be a primitive $p$-th root
of unity. We use that $(1-\zeta)\ZZ[\zeta]\cap\ZZ=p\ZZ$. Since $\ZZ[\zeta]$ is a free $\ZZ$-module, the same holds
coefficientwise for polynomials: if $T\in\ZZ[\mathbf a][X,Z]$ (with any set of variables $\mathbf a$) lies in
$(1-\zeta)\ZZ[\zeta][\mathbf a][X,Z]$, then $T\in p\ZZ[\mathbf a][X,Z]$.

Let $\mathbf a=(a_0,\dots,a_D)$ be variables and $A=\sum a_iX^iZ^{D-i}$ the generic form. The product
$\prod_{j=0}^{p-1}A(\zeta^jX,Z)$ is invariant under $\mathrm{Gal}(\QQ(\zeta)/\QQ)$, which permutes the factors, so its
coefficients lie in $\ZZ[\mathbf a]$. It is also invariant under $X\mapsto\zeta X$, so only monomials $X^iZ^{pD-i}$ with
$p\mid i$ occur. Hence there is a unique form $A^{[p]}$ of degree $D$ over $\ZZ[\mathbf a]$ with
$$A^{[p]}(X^p,Z^p)=\prod_{j=0}^{p-1}A(\zeta^jX,Z).$$
For a form over any ring we define $A^{[p]}$ by specialising the coefficients. Over a field, from
$\prod_j(\zeta^ja-b)=(-1)^{p+1}(a^p-b^p)$,
\begin{equation}\label{eq:frobform}
\Bigl(\prod_{i=1}^D(z_iX-x_iZ)\Bigr)^{[p]}=\bigl((-1)^{p+1}\bigr)^D\prod_{i=1}^D\bigl(z_i^pX-x_i^pZ\bigr).
\end{equation}
We also write $A^{(p)}=\sum a_i^pX^iZ^{D-i}$.

The form $A^{[p]}$ has as zeros the $p$-th powers of the zeros of $A$. The next lemma shows that it is congruent to $A^{(p)}$ modulo $p$, and that it raises congruences by one power of $I$.

\begin{lemma}\label{lem:frob}
Let $A,B$ be forms of degree $D$ over $R$.
\begin{enumerate}[label=\rm(\alph*)]
\item $A^{[p]}\equiv A^{(p)}\pmod{pR}$.
\item If $m\ge1$ and $A\equiv B\pmod{I^m}$, then $A^{[p]}\equiv B^{[p]}\pmod{I^{m+1}}$.
\end{enumerate}
\end{lemma}

\begin{proof}
(a) It suffices to treat generic coefficients. Modulo $1-\zeta$ every factor $A(\zeta^jX,Z)$ is congruent to $A(X,Z)$,
and $A(X,Z)^p\equiv A^{(p)}(X^p,Z^p)$ modulo $p$. Hence $A^{[p]}(X^p,Z^p)-A^{(p)}(X^p,Z^p)$, which has coefficients in
$\ZZ[\mathbf a]$, lies in $(1-\zeta)\ZZ[\zeta][\mathbf a][X,Z]$, so in $p\ZZ[\mathbf a][X,Z]$.

(b) Let $H=\sum y_iX^iZ^{D-i}$ be a second generic form, and decompose
$$\prod_{j=0}^{p-1}(A+H)(\zeta^jX,Z)=\Pi_0+\Pi_1+\Pi_{\ge2},$$
where $\Pi_1$ collects the terms of degree $1$ in $\mathbf y$ and $\Pi_{\ge2}$ those of degree $\ge2$. Galois conjugation
preserves the degree in $\mathbf y$, so each part has coefficients in $\ZZ[\mathbf a,\mathbf y]$. Moreover
$\Pi_1=\sum_jH(\zeta^jX,Z)\prod_{l\ne j}A(\zeta^lX,Z)\equiv pHA^{p-1}\equiv0$ modulo $1-\zeta$, so
$\Pi_1\in p\ZZ[\mathbf a,\mathbf y][X,Z]$. Now specialise so that $A+H=B$; the coefficients of $H$ lie in $I^m$. Then $\Pi_1$
is linear in them with coefficients in $pR\subseteq I$, so its coefficients lie in $I^{m+1}$; and every coefficient of
$\Pi_{\ge2}$ lies in $I^{2m}\subseteq I^{m+1}$, as $m\ge1$. Hence $B^{[p]}(X^p,Z^p)-A^{[p]}(X^p,Z^p)=\Pi_1+\Pi_{\ge2}$ has
coefficients in $I^{m+1}$.
\end{proof}

Next we compare the push-forwards by two pairs of forms that are congruent modulo $I$.

\begin{lemma}\label{lem:pushmod}
Let $\Phi,\Psi$ be pairs of forms of degree $\delta$ over $R$ with $\Phi\equiv\Psi\pmod I$, and $A,B$ forms of degree $D$ with
$A\equiv B\pmod{I^m}$, $m\ge0$. Then $\Phi_*A-\Psi_*A$ has coefficients in $I$, and
$$(\Phi_*A-\Psi_*A)-(\Phi_*B-\Psi_*B)\ \text{ has coefficients in }I^{m+1}.$$
\end{lemma}

\begin{proof}
Write $\Phi=\Psi+H$, with the coefficients $\mathbf h$ of $H$ in $I$. The coefficients of $\Phi_*A$ are integer polynomials
$\mathcal R(\mathbf a,\boldsymbol\psi+\mathbf h)$, and by Taylor expansion
$\mathcal R(\mathbf a,\boldsymbol\psi+\mathbf h)-\mathcal R(\mathbf a,\boldsymbol\psi)=\sum_{\beta\ne0}\mathbf h^\beta\mathcal S_\beta(\mathbf a,\boldsymbol\psi)$
with integer polynomials $\mathcal S_\beta$. Each $\mathbf h^\beta$ lies in $I$, which gives the first claim. As
$\mathbf b-\mathbf a$ has entries in $I^m$, $\mathcal S_\beta(\mathbf b,\boldsymbol\psi)-\mathcal S_\beta(\mathbf a,\boldsymbol\psi)\in I^m$, which gives
the second.
\end{proof}

Combining the two previous lemmas, we show that push-forward by $\Phi$ raises the modulus of a congruence by one power of $\p$.

\begin{lemma}\label{lem:raise}
Let $\Phi$, $\Psi$, $c$, $s$ be as in \eqref{eq:H}, and let $A,B\in\OO[X,Z]$ have degree $D$.
\begin{enumerate}[label=\rm(\alph*)]
\item $\Phi_*A\equiv\pm A^{(p^2)}(cX-sZ,\,Z)\pmod\p$, where $A^{(p^2)}=(A^{(p)})^{(p)}$ and the sign depends only on $p$ and
$D$.
\item If $m\ge1$, $\mu\in\OO^\times$ and $A\equiv\mu B\pmod{\p^m}$, then $\Phi_*A\equiv\mu^d\,\Phi_*B\pmod{\p^{m+1}}$.
\end{enumerate}
\end{lemma}

\begin{proof}
First we show the identity
\begin{equation}\label{eq:psipush}
\Psi_*A=\varepsilon\,(A^{[p]})^{[p]}(cX-sZ,\,Z),\qquad\varepsilon=\pm1\text{ depending only on }p,D,
\end{equation}
for every form $A$ over every ring. Both sides are integer polynomials in $(\mathbf a,c,s)$, so it suffices to check it
for complex values. Write $A=\prod_i(z_iX-x_iZ)$. By \eqref{eq:pushprod},
$$\Psi_*A=\pm\prod_i\bigl(cz_i^dX-(x_i^d+sz_i^d)Z\bigr)=\pm\prod_i\bigl(z_i^d(cX-sZ)-x_i^dZ\bigr),$$
and by \eqref{eq:frobform} applied twice, $(A^{[p]})^{[p]}=\pm\prod_i(z_i^dX-x_i^dZ)$. This proves \eqref{eq:psipush}.

We apply the previous lemmas with $R=\OO$ and $I=\p$, which contains $p$.

(a) By Lemma \ref{lem:pushmod}, $\Phi_*A\equiv\Psi_*A\pmod\p$. By Lemma \ref{lem:frob}(a),
$(A^{[p]})^{[p]}\equiv(A^{[p]})^{(p)}\pmod p$, and since $a\equiv b\pmod p$ implies $a^p\equiv b^p\pmod p$, Lemma
\ref{lem:frob}(a) again gives $(A^{[p]})^{(p)}\equiv(A^{(p)})^{(p)}=A^{(p^2)}\pmod p$. As $p\in\p$, substituting into
\eqref{eq:psipush} gives (a).

(b) Put $A'=\mu B$, so $\Phi_*A'=\mu^d\Phi_*B$ and $A\equiv A'\pmod{\p^m}$. Write
$$\Phi_*A-\Phi_*A'=(\Psi_*A-\Psi_*A')+\bigl[(\Phi_*A-\Psi_*A)-(\Phi_*A'-\Psi_*A')\bigr].$$
The bracket has coefficients in $\p^{m+1}$ by Lemma \ref{lem:pushmod}. By Lemma \ref{lem:frob}(b) twice,
$(A^{[p]})^{[p]}\equiv(A'^{[p]})^{[p]}\pmod{\p^{m+2}}$, and the substitution $(X,Z)\mapsto(cX-sZ,Z)$ has coefficients in
$\OO$, so by \eqref{eq:psipush} the first term has coefficients in $\p^{m+2}$.
\end{proof}

\section{The growing congruence}\label{sec:cong-grow}

In this section we show that  along the orbit of $\varphi$ the characteristic forms satisfy a congruence whose modulus grows by one
power of $\p$ at each step.

\subsection{The period}
Let $\sigma_p$ be the map $y\mapsto y^p$ on $k$, and let $\bar c,\bar s\in k$ be the reductions of $c,s$. Define
$\Theta\colon k^\times\times k\to k^\times\times k$ by
$$\Theta(x,y)=\bigl(\sigma_p^2(x)\,\bar c,\ \sigma_p^2(x)\,\bar s+\sigma_p^2(y)\bigr).$$
For a form $B$ over $k$ put $T(B)(X,Z)=B^{\sigma_p^2}(\bar cX-\bar sZ,\,Z)$, where $B^{\sigma_p^2}$ applies $\sigma_p^2$ to the
coefficients. By Lemma \ref{lem:raise}(a), $T$ describes the action of $\Phi_*$ on reductions.

Some power of $T$, and hence of the reduction of $\Phi_*$, is the identity.

\begin{lemma}\label{lem:period}
There is a least integer $q\ge1$ with $\Theta^q(1,0)=(1,0)$ and $f\mid2q$. It depends only on $\bar c$, $\bar s$ and $k$, and
$T^q(B)=B$ for every form $B$ over $k$.
\end{lemma}

\begin{proof}
The map $(x,y,j)\mapsto(\Theta(x,y),j+2)$ of the finite set $k^\times\times k\times\ZZ/f\ZZ$ is a bijection; its inverse is
explicit. So the orbit of $(1,0,0)$ is purely periodic, and its period is the least $q$ as in the statement. By
induction, $T^n(B)=B^{\sigma_p^{2n}}(\gamma_nX-\delta_nZ,\,Z)$ with $(\gamma_n,\delta_n)=\Theta^n(1,0)$: indeed
$$T\bigl(B^{\sigma_p^{2n}}(\gamma_nX-\delta_nZ,Z)\bigr)=B^{\sigma_p^{2n+2}}\bigl(\sigma_p^2(\gamma_n)(\bar cX-\bar sZ)-\sigma_p^2(\delta_n)Z,\ Z\bigr).$$
For $n=q$ this is $B^{\sigma_p^{2q}}$, which equals $B$ because $\sigma_p^f$ is the identity on $k$ and $f\mid2q$.
\end{proof}

\subsection{The congruence}
For a form $A$ over $\OO$ we write $\Phi_*^nA$ for the $n$-fold push-forward (Definition \ref{def:push}).

We now prove the growing congruence along the orbit.

\begin{theorem}\label{thm:cong}
Let $A\in\OO[X,Z]$ be a primitive form, and put $A_n=\Phi_*^nA$. Then:
\begin{enumerate}[label=\rm(\roman*)]
\item every $A_n$ is primitive;
\item for every $n\ge0$ there is $\mu_n\in\{\pm1\}$ with $A_{n+q}\equiv\mu_n\,A_n\pmod{\p^{n+1}}$.
\end{enumerate}
\end{theorem}

\begin{proof}
Let a bar denote reduction modulo $\p$. Lemma \ref{lem:raise}(a) gives $\overline{\Phi_*B}=\pm T(\overline B)$ for every form
$B$ over $\OO$. The map $T$ is injective, as $\sigma_p$ is injective and the substitution is invertible, and
$T(-B)=-T(B)$.

(i) By induction, $\overline{A_{n+1}}=\pm T(\overline{A_n})\ne0$.

(ii) For $n=0$, Lemma \ref{lem:period} gives $\overline{A_q}=\pm T^q(\overline A)=\pm\overline A$. Suppose
$A_{n+q}\equiv\mu_nA_n\pmod{\p^{n+1}}$. Lemma \ref{lem:raise}(b), with $m=n+1$, gives
$A_{n+1+q}\equiv\mu_n^d\,A_{n+1}\pmod{\p^{n+2}}$, and $\mu_n^d=\pm1$.
\end{proof}

\subsection{Minimal forms and resultants}
For $\beta\in\mathbb P^1(\Kbar)$ let $F_\beta\in\OO[X,Z]$ be a primitive form, irreducible over $K$, vanishing at $\beta$; it is unique up
to $\OO^\times$. Its degree is $[K(\beta):K]$, and its zeros are the $K$-conjugates of $\beta$, each simple.

The growing congruence gives a lower bound for the $\p$-adic valuation of the resultant of the minimal forms of two points of the orbit.

\begin{proposition}\label{prop:res}
Let $\alpha\in\mathbb P^1(\Kbar)$ be not preperiodic for $\varphi$, and $F_n=F_{\alpha_n}$. Let $0\le n<m$ with $m\equiv n\pmod q$. Then
$\Res(F_m,F_n)\ne0$, and
$$v_\p\bigl(\Res(F_m,F_n)\bigr)\ \ge\ (n+1)\,\frac{D_nD_m}{D_0},$$
where $v_\p$ is normalised by $v_\p(\p)=1$.
\end{proposition}

\begin{proof}
Put $A=F_0$ and $A_n=\Phi_*^nA$. We first show $A_n=\kappa_nF_n^{D_0/D_n}$ with $\kappa_n\in\OO^\times$. By \eqref{eq:pushprod}, as
$\Phi_0,\Phi_1$ have no common zero, the zeros of $A_n$ are the points $\varphi^n(\gamma)$, where $\gamma$ runs over the zeros of $A$,
that is, over the $K$-conjugates of $\alpha$. As $K(\alpha_n)\subseteq K(\alpha)$, each $K$-conjugate of $\alpha_n$ occurs $D_0/D_n$
times. So $A_n$ and $F_n^{D_0/D_n}$ are forms over $K$ of the same degree with the same zeros, and $A_n=\kappa_nF_n^{D_0/D_n}$ with
$\kappa_n\in K^\times$. Both are primitive, by Theorem \ref{thm:cong}(i) and Gauss's lemma, so $\kappa_n\in\OO^\times$.

Write $m=n+jq$. Applying Theorem \ref{thm:cong}(ii) at $n,n+q,\dots,n+(j-1)q$, and using $\p^{n+iq+1}\subseteq\p^{n+1}$, we get
$A_m\equiv\pm A_n\pmod{\p^{n+1}}$. By Lemma \ref{lem:resmod}, $v_\p(\Res(A_n,A_m))\ge(n+1)D_0$. By multiplicativity,
$\Res(A_n,A_m)$ equals a unit times $\Res(F_n,F_m)^{(D_0/D_n)(D_0/D_m)}$, which gives the inequality.

If $\Res(F_m,F_n)=0$, then $\alpha_m=\gamma(\alpha_n)$ for some $\gamma\in\mathrm{Gal}(\Kbar/K)$. Since $\hh_\varphi\circ\varphi=d\,\hh_\varphi$ and
$\hh_\varphi$ is invariant under $\mathrm{Gal}(\Kbar/K)$, this gives $d^m\hh_\varphi(\alpha)=d^n\hh_\varphi(\alpha)$, so $\hh_\varphi(\alpha)=0$ and
$\alpha$ is preperiodic, a contradiction.
\end{proof}

\section{Green functions and the proof of Theorem \ref{thm:window}}\label{sec:green}

In this section we compare  the congruence of \S\ref{sec:cong-grow} with an archimedean energy estimate over a range of degrees of iterates of
$\asymp\sqrt{\log D}$ steps.

\subsection{Absolute values and Green functions}\label{sec:greenfn}
For each place $\ell\le\infty$ of $\QQ$ let $\CC_\ell$ be $\CC$ if $\ell=\infty$, and the completion of an algebraic closure of
$\QQ_\ell$ otherwise. For a number field $L'$, every embedding $\sigma\colon L'\hookrightarrow\CC_\ell$ gives an absolute value
$|x|_\sigma=|\sigma x|_\ell$, where $|\ell|_\ell=\ell^{-1}$ and $|\cdot|_\infty$ is the usual absolute value, and a place $w$ of $L'$ above $\ell$ arises from exactly $[L'_w:\QQ_\ell]$ embeddings. We write
$\sum_\sigma$ for the sum over all embeddings of $L'$ into all the fields $\CC_\ell$. The product formula reads
$\sum_\sigma\log|x|_\sigma=0$ for $x\in L'^\times$, and the height of $\beta\in\mathbb P^1(L')$ with a lift $\tilde v\in L'^2$ is
$h(\beta)=[L':\QQ]^{-1}\sum_\sigma\log\|\sigma\tilde v\|$, with $\|\cdot\|$ the max norm.

For an embedding $\iota\colon K\hookrightarrow\CC_\ell$ let $\Phi^\iota$ be the pair obtained by applying $\iota$ to the
coefficients, and for $v\in\CC_\ell^2\setminus\{0\}$ put $\widehat H_\iota(v)=\lim_nd^{-n}\log\|(\Phi^\iota)^{(n)}(v)\|$ and
$r_\iota=\frac1{d(d-1)}\log|\iota\Res(\Phi)|$. Following Baker \cite[(1.6)]{Bak}, for distinct $z,z'\in\mathbb P^1(\CC_\ell)$ with
lifts $\tilde z,\tilde z'$ put
$$g_\iota(z,z')=-\log|\tilde z\wedge\tilde z'|+\widehat H_\iota(\tilde z)+\widehat H_\iota(\tilde z')-r_\iota,$$
where $(v_0,v_1)\wedge(v_0',v_1')=v_0v_1'-v_1v_0'$; by (G1) below it does not depend on the lifts. Let $\mathcal S$ be the finite set of places of
$K$ that are archimedean, or at which some coefficient of $\Phi$ is not integral, or which divide $\Res(\Phi)$. By
\S\ref{sec:setting}, $\p\notin\mathcal S$. We use the following facts.
\begin{enumerate}[label=(G\arabic*),leftmargin=*]
\item $\widehat H_\iota(\Phi^\iota(v))=d\,\widehat H_\iota(v)$ and $\widehat H_\iota(\lambda v)=\log|\lambda|+\widehat H_\iota(v)$.
\item If $\iota$ induces a place of $K$ outside $\mathcal S$, then $\widehat H_\iota(v)=\log\|v\|$, $r_\iota=0$ and
$g_\iota(z,z')=-\log\bigl(|\tilde z\wedge\tilde z'|/\|\tilde z\|\|\tilde z'\|\bigr)\ge0$.
\item For a number field $L'\supseteq K$, $\beta\in\mathbb P^1(L')$ and a lift $\tilde v\in L'^2$,
$\sum_\sigma\widehat H_{\sigma|_K}(\sigma\tilde v)=[L':\QQ]\,\hh_\varphi(\beta)$, where $\hh_\varphi(\beta)=\lim_nd^{-n}h(\varphi^n(\beta))$.
\item (Baker \cite[Thm.~1.1]{Bak}) For every $\iota$ there is $C_\iota\ge0$, with $C_\iota=0$ when $\iota$ induces a place outside
$\mathcal S$, such that $\sum_{i\ne j}g_\iota(z_i,z_j)\ge-C_\iota M\log M$ for all $M\ge2$ distinct $z_1,\dots,z_M\in\mathbb P^1(\CC_\ell)$.
\end{enumerate}
(G1)--(G3) are standard properties of the Call--Silverman canonical height and its local decomposition \cite{CS}; in (G2)
the inequality is the ultrametric inequality, and in (G3) the limit may be exchanged with the sum because, by (G2), only
finitely many places contribute anything other than $\log\|\sigma\tilde v\|$. Baker's Theorem 1.1 holds over any algebraically closed field complete with respect to a non-trivial absolute value,
in particular over each $\CC_\ell$, with a constant depending on the map and the field, and its sum is over ordered pairs; for $\iota$
outside $\mathcal S$, (G4) follows from (G2). We call $\sum_{i\ne j}g_\iota(z_i,z_j)$ the \emph{energy} of $z_1,\dots,z_M$ at $\iota$.
Put
$$C=\max\Bigl\{1,\ \frac1{[K:\QQ]}\sum_\iota C_\iota\Bigr\},\qquad\lambda_\p=\frac{\log N\p}{[K:\QQ]},$$
the sum running over all embeddings of $K$ into all $\CC_\ell$; only finitely many terms are non-zero. We fix the constants $C_\iota$ once and for all; all constants below depend on this choice.

\subsection{Choice of parameters}
Put $B=\sqrt{12qC/\lambda_\p}$, so that $\lambda_\p B^2/(3q)=4C$, let $k_1\ge1$ be the least integer with
$\lambda_\p k_1^2/(3q)\ge C\log(3q+k_1+2+B)$, and define
$$k(D)=3q+k_1+\bigl\lceil B\sqrt{\log3D}\,\bigr\rceil,\qquad\mathcal B(D)=\frac{c_1'}D\exp\bigl(-c_2\sqrt{\log3D}\bigr),$$
with $c_1'=\frac34d^{-(3q+k_1+1)}$ and $c_2=B\log d$.

In the next lemma we prove the inequality that dictates the choice of the length $k(D)$.

\begin{lemma}\label{lem:k}
For every integer $D\ge1$ and $k=k(D)$,
$$\frac{\lambda_\p}{3q}(k+2)(k-3q)-C\log\bigl((k+1)D\bigr)\ \ge\ 2C\log3D.$$
\end{lemma}

\begin{proof}
Put $u=k-3q\ge k_1+B\sqrt{\log3D}$. Then $(k+2)(k-3q)\ge u^2\ge k_1^2+B^2\log3D$, so the first term is at least
$\lambda_\p k_1^2/(3q)+4C\log3D$. Next, $k+1\le3q+k_1+2+B\sqrt{\log3D}\le(3q+k_1+2+B)(1+\log3D)$, and $\log(1+x)\le x$, so
$\log((k+1)D)\le\log D+\log(3q+k_1+2+B)+\log3D\le2\log3D+\log(3q+k_1+2+B)$. The choice of $k_1$ gives the claim.
\end{proof}

\subsection{A bound over iterates of constant degree}

We first prove the  bound when there is is a range of iterates of constant degree that starts at the first step.

\begin{proposition}\label{prop:window}
Let $\alpha\in\mathbb P^1(\Kbar)$ be not preperiodic for $\varphi$, and suppose that $D_n=D_0$ for $0\le n\le k(D_0)$. Then
$\hh_\varphi(\alpha)\ge\mathcal B(D_0)$.
\end{proposition}

\begin{proof}
Write $D=D_0$ and $k=k(D)$. Let $S_n$ be the set of the $D$ distinct $K$-conjugates of $\alpha_n$, the zeros of $F_n=F_{\alpha_n}$,
and $S=S_0\cup\dots\cup S_k$. Let $L'$ be a finite Galois extension of $K$ over which all points of $S$ are defined, and for
each $\beta\in S$ fix a lift $\tilde v_\beta\in L'^2$. For an embedding $\sigma$ of $L'$ write
$g_\sigma(\beta,\beta')=g_{\sigma|_K}(\sigma\beta,\sigma\beta')$.

A conjugate of $\alpha_n$ is never a conjugate of $\alpha_m$ for $n\ne m$ (end of the proof of
Proposition \ref{prop:res}). So $S$ has $N=(k+1)D$ elements, and $N\ge2$ since $k\ge1$.

We claim
\begin{equation}\label{eq:energy}
\sum_\sigma\sum_{\beta\ne\beta'\in S}g_\sigma(\beta,\beta')=2(N-1)\,[L':\QQ]\,D\sum_{n=0}^k\hh_\varphi(\alpha_n).
\end{equation}
Indeed, the terms $-\log|\tilde v_\beta\wedge\tilde v_{\beta'}|_\sigma$ sum to $0$ over $\sigma$ by the product formula, as
$\tilde v_\beta\wedge\tilde v_{\beta'}\in L'^\times$, and so do the terms $r_{\sigma|_K}$, as $\Res(\Phi)\in K^\times$. Each $\beta$ occurs in
$2(N-1)$ ordered pairs, so by (G3) the terms $\widehat H$ contribute $2(N-1)[L':\QQ]\sum_{\beta\in S}\hh_\varphi(\beta)$, and
$\hh_\varphi(\beta)=\hh_\varphi(\alpha_n)$ for $\beta\in S_n$.

Each embedding $\iota$ of $K$ extends to $[L':K]$ embeddings $\sigma$ of $L'$, and each
$\sigma(S)$ consists of $N$ distinct points. By (G2) and (G4), the embeddings whose restriction to $K$ induces a place in
$\mathcal S$ contribute at least $-[L':K]\sum_\iota C_\iota N\log N\ge-[L':\QQ]\,C\,N\log N$ to the left side of \eqref{eq:energy}, and
for all other embeddings every pair contributes $\ge0$.

Let $\sigma$ be an embedding of $L'$ into $\CC_p$ inducing $\p$ on $K$, and let $0\le n<m\le k$
with $m\equiv n\pmod q$. The points of $\sigma(S_m)$ and $\sigma(S_n)$ are the zeros of $\sigma F_m$ and $\sigma F_n$, two forms with
coefficients in the valuation ring of $\CC_p$ and with a unit coefficient. Since $g_\sigma$ does not depend on the lifts, we
may compute it with lifts coming from factorisations of these forms over $\CC_p$ into linear factors $z_iX-x_iZ$, which may be
chosen with $\|(x_i,z_i)\|=1$, because the Gauss norm on $\CC_p[X,Z]$ is multiplicative. Then (G2) and \eqref{eq:resprod} show that the pairs with one point at
level $m$ and the other at level $n$, in both orders, contribute
$$-2\log\bigl|\Res(F_m,F_n)\bigr|_\sigma=\frac{2\log p}e\,v_\p\bigl(\Res(F_m,F_n)\bigr)\ \ge\ \frac{2\log p}e\,(n+1)D$$
by Proposition \ref{prop:res} with $D_n=D_m=D$. All other pairs contribute $\ge0$. The number of embeddings of $L'$ into
$\CC_p$ inducing $\p$ on $K$ is $[L':K]\,[K_\p:\QQ_p]=[L':K]\,ef$. So these embeddings contribute at least
$$2[L':K]\,D\log N\p\sum_{\substack{0\le n<m\le k\\ m\equiv n\ (q)}}(n+1)\ \ge\ 2[L':K]\,D\log N\p\cdot\frac{(k+1)(k+2)(k-3q)}{6q},$$
using $\sum_{n<m\le k,\ q\mid m-n}(n+1)\ge\sum_{n=0}^k(n+1)\bigl(\frac{k-n}q-1\bigr)=\frac{(k+1)(k+2)(k-3q)}{6q}$.

On the right of \eqref{eq:energy}, $\hh_\varphi(\alpha_n)=d^n\hh_\varphi(\alpha)$, so
$\sum_{n\le k}\hh_\varphi(\alpha_n)\le\frac43d^k\hh_\varphi(\alpha)$, as $d\ge4$. Combining, dividing by $[L':\QQ](k+1)D$ and using
$N-1<(k+1)D$,
$$\tfrac83\,D\,d^k\,\hh_\varphi(\alpha)\ \ge\ \frac{\lambda_\p}{3q}(k+2)(k-3q)-C\log\bigl((k+1)D\bigr)\ \ge\ 2C\log3D\ \ge\ 2$$
by Lemma \ref{lem:k}, as $C\ge1$ and $\log3>1$. Since $d^k\le d^{3q+k_1+1}\exp(B\log d\sqrt{\log3D})$, this gives
$\hh_\varphi(\alpha)\ge\mathcal B(D)$.
\end{proof}

\begin{proof}[Proof of Theorem \ref{thm:window}]
Take $k_0=3q+k_1$. Proposition \ref{prop:window} applies to $\alpha_m$, whose orbit is $(\alpha_{m+i})_{i\ge0}$, and gives
$\hh_\varphi(\alpha_m)\ge\mathcal B(D_m)$. As $D_m\le D_0$, $\mathcal B(D_m)\ge\mathcal B(D_0)D_0/D_m$. Now $\hh_\varphi(\alpha_m)=d^m\hh_\varphi(\alpha)$ and
$D_0/D_m=\prod_{n<m}r_n$.
\end{proof}

\section{The Latt\`es map of multiplication by \texorpdfstring{$p$}{p}}\label{sec:redp}

In this section we show that the map induced by $[p]$ on the $x$-line satisfies the hypothesis \eqref{eq:H}, and that it has no periodic critical point. From now on $E/K$ is as in the introduction. 

\subsection{Division polynomials}

Let $E$ be given by $y^2+a_1xy+a_3y=x^3+a_2x^2+a_4x+a_6$ with $a_i\in\OO$. We use the division polynomials $\psi_n$ and
$\phi_n=x\psi_n^2-\psi_{n+1}\psi_{n-1}$ as in \cite[Ex.~3.7]{SilAEC}. We need the following facts. The polynomials
$\psi_n^2$ and $\phi_n$ lie in $\ZZ[a_1,\dots,a_6][x]$; $\phi_n$ is monic of degree $n^2$, and $\psi_n^2$ has degree at
most $n^2-1$ with $x^{n^2-1}$-coefficient $n^2$ \cite[Ex.~3.7]{SilAEC}. Their resultant, taken with these formal degrees,
is
\begin{equation}\label{eq:resdiv}
\Res(\phi_n,\psi_n^2)=\pm\Delta^{n^2(n^2-1)/6}
\end{equation}
(see \cite{CH} and \cite[Ex.~3.7]{SilAEC}). The identity
$x\circ[n]=\phi_n/\psi_n^2$ holds on the universal elliptic curve over $\ZZ[a_1,\dots,a_6,\Delta^{-1}]$. By
\eqref{eq:resdiv}, after specialising to a Weierstrass equation over any field with $\Delta\ne0$, the polynomials $\phi_n$
and $\psi_n^2$ still have no common zero, so the identity specialises: over every field, for every Weierstrass equation,
$\psi_n^2$ is not the zero polynomial and $x\circ[n]=\phi_n/\psi_n^2$ as rational functions on the curve.

Recall that $d=p^2$, and put
$$\Phi_0(X,Z)=Z^d\phi_p(X/Z),\qquad\Phi_1(X,Z)=Z^d\psi_p^2(X/Z).$$
These are binary forms of degree $d$, and $Z$ divides $\Phi_1$.

In characteristic $0$ the forms $\Phi_0,\Phi_1$ define the Latt\`es map of $[p]$.

\begin{lemma}\label{lem:latt}
Over a field of characteristic $0$ the forms $\Phi_0,\Phi_1$ have no common zero in $\mathbb P^1$, and
$x([p]P)=\Phi(x(P))$ for every point $P$, where $x(O)=(1:0)$.
\end{lemma}

\begin{proof}
The function $x\circ[p]$ has degree $2p^2$ on $E$, because $[p]$ has degree $p^2$ and $x$ has degree $2$. It equals $g\circ
x$ with $g=\phi_p/\psi_p^2$, so $g$ has degree $p^2$. If $\phi_p$ and $\psi_p^2$ had a common factor of degree $m\ge1$, then
after cancelling it the numerator would have degree $p^2-m$ and the denominator degree at most $p^2-1-m$, so $g$ would
have degree $p^2-m<p^2$. Hence $\phi_p$ and $\psi_p^2$ are coprime, and since $\Phi_0(1,0)=1$ the forms have no common zero.
The formula for $x([p]P)$ holds for $P\ne O$ with $\psi_p^2(x(P))\ne0$; for the other points both sides are $(1:0)$, using
$\Phi(1,0)=(1,0)$ and, for $P\in E[p]\setminus\{O\}$, that $\phi_p(x(P))\ne0$ by coprimality.
\end{proof}

On a supersingular curve in characteristic $p$, multiplication by $p$ is purely inseparable, which forces the following shape of the map it induces on the $x$-line.

\begin{lemma}\label{lem:redp}
Let $k$ be a field of characteristic $p$ and $\widetilde E/k$ a supersingular elliptic curve given by a Weierstrass
equation. Then there are $c\in k^\times$ and $s\in k$ with
$$\phi_p=x^{p^2}+s,\qquad\psi_p^2=c\qquad\text{in }k[x].$$
Equivalently $\Phi_0=X^d+sZ^d$ and $\Phi_1=cZ^d$, and $x\circ[p]=(x^{p^2}+s)/c$.
\end{lemma}

\begin{proof}
Since $\widetilde E$ is supersingular, $\widetilde E[p](\bar k)=0$ \cite[V.3.1]{SilAEC}. The isogeny $[p]$ has degree $p^2$
\cite[III.6.2]{SilAEC}, and its separable degree equals the number of points of its kernel \cite[III.4.10]{SilAEC}, which is
$1$. So $[p]$ is purely inseparable of degree $p^2$, and by \cite[II.2.12]{SilAEC} it factors as $[p]=\lambda\circ\pi$,
where $\pi\colon\widetilde E\to\widetilde E^{(p^2)}$ is the $p^2$-power Frobenius and $\lambda$ is an isogeny of degree $1$,
hence an isomorphism of elliptic curves. An isomorphism between curves given by Weierstrass equations has the form
$x=u^2x'+r$ with $u\in\bar k^\times$, $r\in\bar k$ \cite[III.3.1]{SilAEC}. Since $x'\circ\pi=x^{p^2}$, we get
$x\circ[p]=u^2x^{p^2}+r$ in $\bar k(x)$. On the other hand $x\circ[p]=\phi_p/\psi_p^2$ with $\psi_p^2\ne0$. Thus
$\phi_p=(u^2x^{p^2}+r)\psi_p^2$. The left side is monic of degree $p^2$, and if $\psi_p^2$ had degree $m\ge1$ the right side
would have degree $p^2+m$. So $\psi_p^2=c\in k^\times$, and comparing coefficients gives $cu^2=1$ and $\phi_p=x^{p^2}+s$ with
$s=cr\in k$.
\end{proof}

\subsection{The Latt\`es map lifts Frobenius}
The Weierstrass equation of $E$ has coefficients in $\OO$, so all coefficients of $\Phi$ lie in $\OO$, and the reduction of
$\Phi$ is the pair attached in the same way to $\widetilde E$. Choose lifts $c,s\in\OO$ of the elements given by Lemma
\ref{lem:redp} for $\widetilde E$. Then $c\in\OO^\times$ and
$$\Phi\ \equiv\ \Psi=\bigl(X^d+sZ^d,\ cZ^d\bigr)\pmod\p.$$
So $\Phi$ satisfies \eqref{eq:H}, and all results of \S\S\ref{sec:setting}--\ref{sec:green} apply to the map $\varphi$ defined by $\Phi$, the \emph{Latt\`es map}
of $[p]$. By Lemma \ref{lem:latt}, $x([p]Q)=\varphi(x(Q))$ for all $Q\in E(\Kbar)$. For $P\in E(\Kbar)$ we write
$$\alpha_n=\varphi^n(x(P))=x([p^n]P),\qquad D_n=[K(\alpha_n):K],\qquad r_n=D_n/D_{n+1},$$
as in \S\ref{sec:setting}.

The canonical height of the Latt\`es map is twice the N\'eron--Tate height.

\begin{lemma}\label{lem:hE}
For every $Q\in E(\Kbar)$, $\hh_\varphi(x(Q))=2\hh_E(Q)$. In particular $x(Q)$ is preperiodic for $\varphi$ if and only if $Q$ is a
torsion point.
\end{lemma}

\begin{proof}
The N\'eron--Tate height satisfies $\hh_E(Q)=\frac12h(x(Q))+O(1)$ and $\hh_E([p]Q)=p^2\hh_E(Q)$ \cite[VIII.9.3]{SilAEC}. Hence
$d^{-n}h(\varphi^n(x(Q)))=p^{-2n}h(x([p^n]Q))=2\hh_E(Q)+O(p^{-2n})$. The last assertion follows, as $\hh_E$ vanishes exactly at the
torsion points. Other normalisations of $\hh_E$ change only the constants in our results, not the exponents.
\end{proof}

Finally, the critical points of the Latt\`es map come from $2$-torsion, and none of them is periodic.

\begin{lemma}\label{lem:crit}
No critical point of the Latt\`es map $\varphi$ is periodic.
\end{lemma}

\begin{proof}
The map $x\colon E\to\mathbb P^1$ has degree $2$ and is ramified exactly at $E[2]$, and $[p]$ is unramified, as we are in characteristic
$0$. For $Q\in E(\Kbar)$, writing $\mathrm{ram}_g(z)$ for the ramification index of a map $g$ at $z$, the identity $\varphi\circ x=x\circ[p]$ gives
$\mathrm{ram}_\varphi(x(Q))\,\mathrm{ram}_x(Q)=\mathrm{ram}_x([p]Q)$. So $x(Q)$ is a
critical point of $\varphi$ exactly when $[p]Q\in E[2]$ and $Q\notin E[2]$. Then $\varphi^n(x(Q))=x([p^n]Q)\in x(E[2])$ for $n\ge1$, and
$x(E[2])$ is stable under $\varphi$; but $x(Q)\notin x(E[2])$, as $x(Q)=x(T)$ forces $Q=\pm T$. Hence $x(Q)$ is not periodic. Finally, every critical point of $\varphi$ in $\mathbb P^1(\CC_p)$ is a zero of the
non-zero form $\partial_X\Phi_0\,\partial_Z\Phi_1-\partial_Z\Phi_0\,\partial_X\Phi_1$ over $K$, hence lies in $\mathbb P^1(\Kbar)$.
\end{proof}

\section{Local dynamics at \texorpdfstring{$\p$}{p}}\label{sec:local}

In this section $\Phi$, $c$, $s$ and $\varphi$ are as in \eqref{eq:H}, for instance the Latt\`es pair and map of
\S\ref{sec:redp}, and $d=p^2$. We describe the dynamics of $\varphi$ on the residue discs of $\mathbb P^1(\CC_p)$: the reduction
of $\varphi$ permutes the residue classes, each residue disc contains exactly one periodic point, and orbits approach these
periodic points at a rate that depends on the ramification.

Let $\CC_p$ be the completion of an algebraic closure of $K_\p$, fix an embedding $\Kbar\hookrightarrow\CC_p$ over $K$, let $v$ be
the valuation with $v(p)=1$, and $\OO_{\CC_p}$ the valuation ring. For a form or a matrix $H$ over $\CC_p$ let $v(H)$ be the least
valuation of its coefficients. Since $v(\pi)=1/e$ for a uniformiser $\pi$ of $\OO$, the hypothesis \eqref{eq:H} says that $v(\Phi-\Psi)\ge1/e$, where
$$\Psi=\bigl(X^d+sZ^d,\ cZ^d\bigr)=G_0\circ(X^d,Z^d),\qquad G_0=\begin{pmatrix}1&s\\0&c\end{pmatrix}\in\mathrm{GL}_2(\OO).$$
Put
$$\rho(u)=\min\{du,\ u+1/e\}\qquad(u\ge0),$$
an increasing function with $\rho(u)>u$ for $u>0$.

Every $P\in\mathbb P^1(\CC_p)$ has a primitive representative $(x_0,z_0)$, with $\min\{v(x_0),v(z_0)\}=0$, unique up to a unit,
and a reduction $\bar P\in\mathbb P^1(\overline{\FF}_p)$. For $P,Q$ with primitive representatives $(x_0,z_0),(y_0,u_0)$ put
$$w(P,Q)=v(x_0u_0-z_0y_0)\in[0,\infty].$$
Then $p^{-w}$ is the chordal metric, which is an ultrametric; so $w(P,R)\ge\min\{w(P,Q),w(Q,R)\}$. Moreover $w(P,Q)>0$ if and
only if $\bar P=\bar Q$, $w$ is invariant under $\mathrm{GL}_2(\OO_{\CC_p})$, and $\mathbb P^1(\CC_p)$ is complete. For $\xi\in\mathbb P^1(\overline{\FF}_p)$ let $U_\xi=\{P:\bar P=\xi\}$ be its residue
disc.

The first lemma describes $\varphi$ near a point in local coordinates, and introduces the invariant $\vartheta$, which measures how far $\varphi$ is from being ramified at a point.

\begin{lemma}\label{lem:loc}
\begin{enumerate}[label=\rm(\alph*)]
\item The reduction $\bar\varphi$ is a bijection of $\mathbb P^1(\overline{\FF}_p)$, and every point of
$\mathbb P^1(\overline{\FF}_p)$ is periodic for it; let $q_\xi$ be the exact period of $\xi$.
\item For $P\in\mathbb P^1(\CC_p)$ choose $M_1,M_2\in\mathrm{GL}_2(\OO_{\CC_p})$ with $M_1(0)=P$ and $M_2(0)=\varphi(P)$, where
$0=[0:1]$, and let $(U_0,U_1)=M_2^{-1}\circ\Phi\circ M_1$, $a(t)=U_0(t,1)$, $b(t)=U_1(t,1)$. Then $a(0)=0$, $b(t)$ is a unit for
$v(t)>0$, and
$$a(t)=u\,t^d+t\,a_1(t),\qquad u\in\OO_{\CC_p}^\times,\quad a_1\in\OO_{\CC_p}[t],\quad v(a_1)\ge1/e.$$
\item For $P,Q$ with $\bar P=\bar Q$: $w(\varphi P,\varphi Q)\ge\rho(w(P,Q))$.
\item In the notation of (b), $\vartheta(P)=v(a'(0))$ does not depend on the choice of $M_1,M_2$, and
$\vartheta(P)=v(J(x_0,z_0))-2$ for any primitive representative $(x_0,z_0)$ of $P$, where
$J=\partial_X\Phi_0\,\partial_Z\Phi_1-\partial_Z\Phi_0\,\partial_X\Phi_1$. One has $\vartheta(P)\ge1/e$, $\vartheta(P)=\infty$ exactly
when $P$ is a critical point of $\varphi$, and there is $\vartheta_0$ with $\vartheta(P)=\vartheta_0$ for all $P$ outside finitely
many residue discs.
\end{enumerate}
\end{lemma}

\begin{proof}
(a) The reduction is $\bar G_0\circ\mathrm{Frob}_p^2$, a composite of bijections of $\mathbb P^1(\overline{\FF}_p)$ which preserve
$\mathbb P^1(\FF)$ for every finite field $\FF$ containing $k$; on the finite set $\mathbb P^1(\FF)$ a bijection is a permutation.

(b) For $M_1=\begin{psmallmatrix}m_{00}&m_{01}\\m_{10}&m_{11}\end{psmallmatrix}$ one has
$(m_{00}X+m_{01}Z)^d\equiv m_{00}^dX^d+m_{01}^dZ^d$ modulo $p\,\OO_{\CC_p}[X,Z]$, since $p\mid\binom di$ for $0<i<d$. Hence
$\Psi\circ M_1\equiv G_0M_1^{(d)}(X^d,Z^d)$ modulo $p$, where $M_1^{(d)}$ has entries $m_{ij}^d$, and since $v(\Phi-\Psi)\ge1/e$ and
$1/e\le1$,
$$(U_0,U_1)\equiv\mu'\,(X^d,Z^d)\pmod{\{v\ge1/e\}},\qquad\mu'=M_2^{-1}G_0M_1^{(d)},$$
with $\det\mu'$ a unit. As $\Res(\Phi_0,\Phi_1)$ is a unit, $(U_0(\tilde P),U_1(\tilde P))$ is primitive for every primitive $\tilde P$;
since $M_2^{-1}\varphi M_1$ fixes $0$, $a(0)=0$ and $b(0)$ is a unit. Now $a(0)\equiv\mu'_{01}$, so $v(\mu'_{01})\ge1/e$ and
$\mu'_{00},\mu'_{11}$ are units. For $v(t)>0$, $b(t)\equiv\mu'_{10}t^d+\mu'_{11}$ modulo the maximal ideal is a unit. Finally
$a(t)-\mu'_{00}t^d$ has all coefficients of valuation $\ge1/e$ and vanishes at $t=0$; so it equals $t\,a_1(t)$ as claimed, with
$u=\mu'_{00}$.

(c) In (b), $Q=M_1([t:1])$ with $v(t)=w(P,Q)>0$, because $w$ is $\mathrm{GL}_2(\OO_{\CC_p})$-invariant and $M_1^{-1}Q$ lies in the
residue disc of $0$. Then $\varphi(Q)=M_2([a(t):b(t)])$ with $b(t)$ a unit, so
$w(\varphi P,\varphi Q)=w([0:1],[a(t):b(t)])=v(a(t))\ge\min\{dv(t),v(t)+1/e\}$.

(d) Replacing $M_1$ by $M_1\circ T$ with $T\in\mathrm{GL}_2(\OO_{\CC_p})$ fixing $0$ multiplies $a'(0)/b(0)$ by $T'(0)$, a unit;
similarly for $M_2$; and $b(0)$ is a unit. By the chain rule $J_{M_2^{-1}\Phi M_1}=\det(M_2)^{-1}\det(M_1)\,J_\Phi\circ M_1$. For
$(U_0,U_1)$, Euler's identities at $(0,1)$ give $\partial_ZU_0(0,1)=d\,a(0)=0$ and $\partial_ZU_1(0,1)=d\,b(0)$, so
$J_{(U_0,U_1)}(0,1)=d\,a'(0)\,b(0)$. As $M_1(0,1)$ is a primitive representative of $P$ and $v(d)=2$, this gives
$\vartheta(P)=v(J(x_0,z_0))-2$. By (b), $a'(0)=a_1(0)$ has valuation $\ge1/e$, and $a'(0)=0$ exactly when $M_2^{-1}\varphi M_1$,
equivalently $\varphi$ at $P$, is ramified. The form $J$ of degree $2d-2$ is non-zero, $\varphi$ being separable; write
$J=\lambda_JJ_0$ with $\lambda_J\in K_\p^\times$ and $J_0$ having a unit coefficient. The reduction $\bar J_0\ne0$ has at most
$2d-2$ zeros in $\mathbb P^1(\overline{\FF}_p)$, and outside their residue discs $\vartheta=v(\lambda_J)-2$.
\end{proof}

Using the contraction in Lemma \ref{lem:loc}(c), we find exactly one periodic point in each residue disc.

\begin{lemma}\label{lem:periodic}
Each residue disc $U_\xi$ contains exactly one periodic point $\eta_\xi$ of $\varphi$. It has exact period $q_\xi$,
$\varphi(\eta_\xi)=\eta_{\bar\varphi(\xi)}$, and if $P\in U_\xi$ is not periodic then $w(P,\eta_\xi)=w(\varphi^{q_\xi}(P),P)$.
\end{lemma}

\begin{proof}
Put $q'=q_\xi$ and $G=\varphi^{q'}$. Then $G(U_\xi)\subseteq U_\xi$, and by Lemma \ref{lem:loc}(c), $w(G(P),G(Q))>w(P,Q)$ for
$P\ne Q$ in $U_\xi$; so $G$ has at most one fixed point in $U_\xi$. For existence, take $P\in U_\xi$ and put
$\delta_j=w(G^{j+1}(P),G^j(P))$. Then $\delta_0>0$ and, applying Lemma \ref{lem:loc}(c) $q'$ times, $\delta_{j+1}\ge\rho(\delta_j)\ge\delta_j+\min\{(d-1)\delta_0,1/e\}$, so $\delta_j\to\infty$. By the
ultrametric inequality the sequence $G^j(P)$ is Cauchy; its limit $\eta$ lies in the closed set $U_\xi$ and is fixed by $G$,
which is continuous. For $0<j<q'$, $\varphi^j(\eta)$ lies in $U_{\bar\varphi^j(\xi)}\ne U_\xi$, so $\eta$ has exact period $q'$. A
periodic point in $U_\xi$ has period divisible by $q'$, since its period fixes $\xi$ under $\bar\varphi$; so it is a fixed point of
some $G^j$, $j\ge1$, and the argument above, applied to $G^j$, shows that it equals $\eta$. The point $\varphi(\eta_\xi)$ is periodic
and lies in $U_{\bar\varphi(\xi)}$, so it equals $\eta_{\bar\varphi(\xi)}$. If $P\in U_\xi$ is not periodic and $u=w(P,\eta_\xi)$, then
$w(G(P),\eta_\xi)=w(G(P),G(\eta_\xi))>u$, so $w(G(P),P)=u$.
\end{proof}

Consequently, the periodic points of $\varphi$ in $\mathbb P^1(\CC_p)$ are exactly the points $\eta_\xi$, $\xi\in\mathbb P^1(\overline{\FF}_p)$.

Orbits of points defined over a finite extension of $K_\p$ approach the orbits of the periodic points of their discs, at a rate controlled by the ramification index.

\begin{lemma}\label{lem:approach}
Let $F\subset\CC_p$ be a finite extension of $K_\p$ with absolute ramification index $e_F$, and put
$m_0=\lceil\log_d(e_F/e)\rceil\ge0$. Let $P\in\mathbb P^1(F)$ have infinite orbit, and $\eta=\eta_{\bar P}$. Then for all $n\ge0$,
$$w\bigl(\varphi^n(P),\varphi^n(\eta)\bigr)\ \ge\ \min\{d^n/e_F,\ 1/e\},$$
and $w(\varphi^n(P),\varphi^n(\eta))\ge(n-m_0+1)/e$ for $n\ge m_0$.
\end{lemma}

\begin{proof}
As $e$ divides $e_F$, $m_0\ge0$. Let $q'=q_{\bar P}$. Choosing primitive representatives with coordinates in $F$,
$w(\varphi^{q'}(P),P)$ is the valuation of a non-zero element of $F$ (non-zero as $P$ is not periodic), and it is positive; so it
is at least $1/e_F$. By Lemma \ref{lem:periodic}, $g_0=w(P,\eta)\ge1/e_F$. Put $g_n=w(\varphi^n(P),\varphi^n(\eta))$; by Lemma
\ref{lem:loc}(c), $g_{n+1}\ge\rho(g_n)$. As $\rho$ is increasing, induction gives $g_n\ge\min\{d^n/e_F,1/e\}$: if this holds for
$n$, then $g_{n+1}\ge\min\{d\min\{d^n/e_F,1/e\},\ 1/e\}\ge\min\{d^{n+1}/e_F,1/e\}$. Since $d^{m_0}\ge e_F/e$, $g_{m_0}\ge1/e$. If
$g\ge1/e$ then $dg\ge g+1/e$, so $\rho(g)=g+1/e$; hence $g_n\ge(n-m_0+1)/e$ for $n\ge m_0$.
\end{proof}

When no critical point is periodic, $\vartheta$ is bounded on the periodic points, which gives the integer $\nu$ used in \S\ref{sec:clusters}.

\begin{lemma}\label{lem:kappa}
Suppose that no critical point of $\varphi$ is periodic. Then
$$\kappa^*=\sup\{\vartheta(\eta):\ \eta\text{ periodic}\}$$
is finite and at least $1/e$, so $\nu=\lfloor e\kappa^*\rfloor$ is an integer with $\nu\ge1$ and $(\nu+1)/e>\kappa^*$. For the
Latt\`es map, $\kappa^*$ and $\nu$ depend only on $E$, $K$ and $\p$.
\end{lemma}

\begin{proof}
By Lemma \ref{lem:loc}(d), $\vartheta\ge1/e$, and $\vartheta(\eta_\xi)=\vartheta_0<\infty$ outside finitely many $\xi$. For each of the
finitely many remaining $\xi$, $\eta_\xi$ is periodic, hence not critical, so $\vartheta(\eta_\xi)<\infty$. Since the periodic points
are exactly the $\eta_\xi$, $\kappa^*$ is a maximum of finitely many finite numbers.
\end{proof}

By Lemma \ref{lem:crit}, Lemma \ref{lem:kappa} applies to the Latt\`es map.

\section{Ramified clusters and the proof of Theorem \ref{thm:main}}\label{sec:clusters}

In this section we show that a drop of the degree forces many conjugates to be highly ramified above $\p$, and use the resulting $\p$-adic clusters to prove Theorem \ref{thm:main}.

In this section $\varphi$ is as in \S\ref{sec:local}, with no periodic critical point, and $\kappa^*$, $\nu$ are as in Lemma
\ref{lem:kappa}. For $\alpha\in\mathbb P^1(\Kbar)$ we write $\alpha_n=\varphi^n(\alpha)$, $D_n=[K(\alpha_n):K]$ and $r_n=D_n/D_{n+1}$, as in
\S\ref{sec:setting}; the degree \emph{drops} at step $n$ if $r_n\ge2$. For a $K$-conjugate $x\in\mathbb P^1(\Kbar)\subset\mathbb P^1(\CC_p)$ of
$\alpha$ let $K_\p(x)\subset\CC_p$ be the field generated over $K_\p$ by the affine coordinate of $x$ (with $K_\p(\infty)=K_\p$), and
$e_x$ its ramification index over $\QQ_p$. It is the completion of $K(\alpha)$ at a prime above $\p$, so $e_x$ is at most the
largest absolute ramification index of a prime of $K(\alpha)$ above $\p$.

\subsection{Drops force ramification}

We first show how the conjugates of $\alpha$ distribute over those of $\alpha_n$.

\begin{lemma}\label{lem:fibres}
Let $\alpha\in\mathbb P^1(\Kbar)$ and $n\ge0$. The map $x\mapsto\varphi^n(x)$ sends the $K$-conjugates of $\alpha$ onto those of
$\alpha_n$, with fibres of exactly $D_0/D_n$ elements; and $\varphi$ sends the $K$-conjugates of $\alpha_n$ onto those of $\alpha_{n+1}$,
with fibres of exactly $r_n$ elements.
\end{lemma}

\begin{proof}
As $\varphi$ is defined over $K$, the maps are $\mathrm{Gal}(\Kbar/K)$-equivariant; they are onto because the Galois group acts
transitively on the conjugates of $\alpha_n$, and it permutes the fibres transitively, so they all have the same size.
\end{proof}

The next lemma is a divided-difference estimate: if two distinct points near a periodic point $\eta$ have the same image, then $\vartheta(\eta)$ is large.

\begin{lemma}\label{lem:divdiff}
Let $\eta\in\mathbb P^1(\CC_p)$, and let $z\ne z'$ in $U_{\bar\eta}$ satisfy $\varphi(z)=\varphi(z')$. Then
$\vartheta(\eta)\ge\min\{w(z,\eta),w(z',\eta)\}$.
\end{lemma}

\begin{proof}
Choose $M_1,M_2$ as in Lemma \ref{lem:loc}(b) for $P=\eta$, and write $z=M_1([t_1:1])$, $z'=M_1([t_2:1])$; as in the proof of Lemma
\ref{lem:loc}(c), $v(t_1)=w(z,\eta)>0$ and $v(t_2)=w(z',\eta)>0$, and $t_1\ne t_2$. Since $\varphi(z)=M_2([a(t_1):b(t_1)])$ and
likewise for $z'$, we have $a(t_1)b(t_2)=a(t_2)b(t_1)$. The polynomial
$$\Gamma(\tau_1,\tau_2)=\frac{a(\tau_1)b(\tau_2)-a(\tau_2)b(\tau_1)}{\tau_1-\tau_2}$$
has coefficients in $\OO_{\CC_p}$, since $\tau_1^i\tau_2^j-\tau_2^i\tau_1^j$ is divisible by $\tau_1-\tau_2$ in $\ZZ[\tau_1,\tau_2]$.
It satisfies $\Gamma(t_1,t_2)=0$ and $\Gamma(0,0)=a'(0)b(0)-a(0)b'(0)=a'(0)b(0)$, with $b(0)$ a unit. Hence
$\vartheta(\eta)=v(\Gamma(0,0)-\Gamma(t_1,t_2))\ge\min\{v(t_1),v(t_2)\}$.
\end{proof}

Combining the previous lemmas, we show that a drop of the degree at step $n$ forces at least half of the conjugates of $\alpha$ to be highly ramified above $\p$.

\begin{lemma}\label{lem:half}
Let $\alpha$ be not preperiodic, and suppose that the degree drops at step $n$. Then the set
$$H_n\ =\ \bigl\{x:\ x\text{ a $K$-conjugate of $\alpha$ with }e_x>e\,d^{\,n-\nu}\bigr\}$$
has at least $D_0/2$ elements.
\end{lemma}

\begin{proof}
Call a conjugate $z$ of $\alpha_n$ \emph{bad} if some conjugate $x$ of $\alpha$ with $\varphi^n(x)=z$ has $e_x\le e\,d^{\,n-\nu}$. Since
$e$ divides $e_x$, there are no bad elements if $n<\nu$.

We claim that each fibre of $\varphi$ on the conjugates of $\alpha_n$ contains at most one bad element. Suppose that $z\ne z'$ in one
fibre are bad, with $x,x'$ as in the definition. Then $\varphi(z)=\varphi(z')$, so $\bar z=\bar z'$, as $\bar\varphi$ is injective (Lemma
\ref{lem:loc}(a)); let $\eta$ be the periodic point of this residue disc. Since $\bar\varphi^n(\bar x)=\bar z$, Lemma \ref{lem:periodic}
gives $\varphi^n(\eta_{\bar x})=\eta$. Put $m_0=\lceil\log_d(e_x/e)\rceil$. As $e_x/e\le d^{\,n-\nu}$ and $n-\nu$ is an integer,
$0\le m_0\le n-\nu$, and Lemma \ref{lem:approach}, applied to $x$ with $F=K_\p(x)$, gives
$$w(z,\eta)\ \ge\ \frac{n-m_0+1}e\ \ge\ \frac{\nu+1}e\ >\ \kappa^*.$$
Likewise $w(z',\eta)>\kappa^*$. By Lemma \ref{lem:divdiff}, $\vartheta(\eta)>\kappa^*$, contradicting the definition of $\kappa^*$.

By Lemma \ref{lem:fibres} the fibres have $r_n\ge2$ elements, so at least $(1-1/r_n)D_n\ge D_n/2$ conjugates of $\alpha_n$ are not bad,
and the $D_0/D_n$ conjugates of $\alpha$ above each of them lie in $H_n$.
\end{proof}

\subsection{A cluster of conjugates}

A conjugate with ramification index $e_x$ lies in a cluster of $e_x/e$ conjugates whose orbits remain $\p$-adically close.

\begin{lemma}\label{lem:cluster}
Let $\alpha$ be not preperiodic, $x$ a $K$-conjugate of $\alpha$, and $K^{\mathrm{ur}}_x\subseteq K_\p(x)$ the maximal unramified extension of $K_\p$
in $K_\p(x)$. The set $\mathcal C(x)$ of conjugates of $x$ over $K^{\mathrm{ur}}_x$ consists of $e_x/e$ conjugates of $\alpha$, and
$$w\bigl(\varphi^n(x'),\varphi^n(x'')\bigr)\ \ge\ \min\{d^n/e_x,\ 1/e\}\qquad\text{for }x'\ne x''\text{ in }\mathcal C(x)\text{ and }n\ge0.$$
\end{lemma}

\begin{proof}
If $x=\infty$, then $\mathcal C(x)=\{x\}$ and there is nothing to prove; let $x\ne\infty$. The set $\mathcal C(x)$ has $[K_\p(x):K^{\mathrm{ur}}_x]=e_x/e$ elements, all roots of the minimal
polynomial of $\alpha$ over $K$. For $x'\in\mathcal C(x)$ the $K^{\mathrm{ur}}_x$-isomorphism $K^{\mathrm{ur}}_x(x)\to K^{\mathrm{ur}}_x(x')$ is isometric, $K^{\mathrm{ur}}_x$ being complete; so
$K_\p(x')=K^{\mathrm{ur}}_x(x')$ has ramification index $e_x$. Since $K_\p(x)$ and $K^{\mathrm{ur}}_x$ have the same residue field, there is $t\in\mathbb P^1(K^{\mathrm{ur}}_x)$ with
$w(x,t)>0$, and then $w(x',t)=w(x,t)>0$. So $\mathcal C(x)$ lies in one residue disc $U_\xi$. Lemma \ref{lem:approach}, applied to each
$x'\in\mathcal C(x)$ with $F=K_\p(x')$ and the same $\eta=\eta_\xi$, and the ultrametric inequality give the claim.
\end{proof}

\subsection{Energy at the last drop}
We use the notation of \S\ref{sec:green}: the Green functions $g_\iota$, the set $\mathcal S$ of places, the properties (G1)--(G4),
and the constants $C$ and $\lambda_\p$.

We now bound $\hh_\varphi(\alpha)$ from below by the energy at $\p$ of the conjugates of $\alpha_m$, where $m$ is a step at which the degree drops.

\begin{proposition}\label{prop:energy}
Let $\alpha$ be not preperiodic, and suppose that the degree drops at step $m$. Then
$$\hh_\varphi(\alpha)\ \ge\ \frac{\lambda_\p}{4d^{\,\nu}D_0}-\frac{e\lambda_\p+2C\log D_m}{4D_m\,d^{\,m}}.$$
\end{proposition}

\begin{proof}
Let $S$ be the set of the $N=D_m$ conjugates of $\alpha_m$; as the degree drops at step $m$, $N\ge r_m\ge2$. Let $L'$ be a finite
Galois extension of $K$ over which the points of $S$ are defined, and for an embedding $\sigma$ of $L'$ write
$g_\sigma(\beta,\beta')=g_{\sigma|_K}(\sigma\beta,\sigma\beta')$. The proof follows the lines of the proof of Proposition \ref{prop:window}: with the
single level $S$ in place of $S_0\cup\dots\cup S_k$,
$$\sum_\sigma\sum_{\beta\ne\beta'\in S}g_\sigma(\beta,\beta')\ =\ 2(N-1)\,[L':\QQ]\,N\,\hh_\varphi(\alpha_m)\ \le\ 2[L':\QQ]\,N^2d^{\,m}\,\hh_\varphi(\alpha),$$
the embeddings inducing places in $\mathcal S$ contribute at least $-[L':\QQ]\,C\,N\log N$, and all others contribute at least
$0$.

Let $\sigma$ be one of the $[L':K]\,ef$ embeddings of $L'$ into $\CC_p$ inducing $\p$ on $K$. By (G2),
$g_\sigma(\beta,\beta')=w(\sigma\beta,\sigma\beta')\log p$. The restriction of $\sigma$ to $K$ and the fixed embedding of $K$ both
extend to $\QQ_p$-embeddings of $K_\p$ into $\CC_p$, so they differ by a continuous automorphism $\theta$ of $\CC_p$, which preserves
$v$; then $\theta^{-1}\circ\sigma$ and the fixed embedding agree on $K$, so they differ on $L'$ by an element of $\mathrm{Gal}(L'/K)$,
which permutes $S$. Hence $\sum_{\beta\ne\beta'\in S}w(\sigma\beta,\sigma\beta')$ equals its value $\mathcal G$ for the fixed
embedding.

By Lemma \ref{lem:half}, $|H_m|\ge D_0/2$. Let $x\in H_m$ and $y=\varphi^m(x)$. By Lemmas \ref{lem:cluster} and \ref{lem:fibres}, the set
$\varphi^m(\mathcal C(x))\subseteq S$ contains $y$ and at least $(e_x/e)N/D_0$ points, pairwise at chordal valuation
$\ge\delta_x=\min\{d^m/e_x,1/e\}\le1$. As $e_x>e\,d^{\,m-\nu}$, the point $y$ contributes to $\mathcal G$ at least
$$\Bigl(\frac{e_xN}{eD_0}-1\Bigr)\delta_x\ \ge\ \frac N{D_0}\min\Bigl\{\frac{d^{\,m}}e,\ \frac{e_x}{e^2}\Bigr\}-1\ \ge\
\frac{d^{\,m-\nu}N}{eD_0}-1.$$
There are at least $|H_m|N/D_0\ge N/2$ such points $y$, so $\mathcal G\ge\frac N2\bigl(\frac{d^{\,m-\nu}N}{eD_0}-1\bigr)$ (if the bracket
is negative, use $\mathcal G\ge0$). Collecting, and dividing by $[L':\QQ]$, as $[L':K]ef\log p/[L':\QQ]=e\lambda_\p$,
$$2N^2d^{\,m}\,\hh_\varphi(\alpha)\ \ge\ e\lambda_\p\mathcal G-CN\log N\ \ge\ \frac{\lambda_\p d^{\,m-\nu}N^2}{2D_0}-\frac{e\lambda_\p N}2-CN\log N,$$
which is the claim.
\end{proof}

\subsection{Proof of Theorem \ref{thm:main}}

\begin{proof}[Proof of Theorem \ref{thm:main}]
Let $k(D)$, $c_1'$, $c_2$ and $\mathcal B$ be as in Theorem \ref{thm:window} for the Latt\`es map $\varphi$, which has no periodic
critical point by Lemma \ref{lem:crit}. Put $c_4=\lambda_\p/(4d^{\,\nu})$ and $c_5=\max\{1,\,e\lambda_\p/4+C/2\}$, so that
$(e\lambda_\p+2C\log N)/4\le c_5\log3N$ for all $N\ge1$. Let $P$ be non-torsion, $\alpha=x(P)$, which is not preperiodic by Lemma
\ref{lem:hE}, $D=D_0$ and $X_D=\exp(-c_2\sqrt{\log3D})$, so that $\mathcal B(D)=c_1'X_D/D$. Recall $\hh_E(P)=\frac12\hh_\varphi(\alpha)$.

If the degree never drops, Theorem \ref{thm:window} with $m=0$ gives $\hh_\varphi(\alpha)\ge c_1'X_D/D$, which implies both (i) and
(ii). Otherwise, as $(D_n)$ is a non-increasing sequence of positive integers, there is a last step $m$ at which the degree drops,
and Theorem \ref{thm:window}, applied with $m+1$ in place of $m$, gives
\begin{equation}\label{eq:lastdrop}
\hh_\varphi(\alpha)\ \ge\ \frac{c_1'X_D}D\prod_{n\le m}\frac{r_n}d\ =\ \frac{c_1'X_D}{D_{m+1}\,d^{\,m+1}}.
\end{equation}

For (i): by Lemma \ref{lem:half} there is a conjugate $x$ of $\alpha$ with $e\,d^{\,m-\nu}<e_x$, and $e_x\le e_\p(P)$ since $K(\alpha)\subseteq
K(P)$; so $d^{\,m+1}<d^{\,\nu+1}e_\p(P)/e\le d^{\,\nu+1}e_\p(P)$, and \eqref{eq:lastdrop} with $D_{m+1}\le D$ gives
$\hh_\varphi(\alpha)\ge c_1'd^{-\nu-1}X_D/(e_\p(P)D)$.

For (ii): by \eqref{eq:lastdrop} and $D_{m+1}\le D_m$, we have $\hh_\varphi(\alpha)\ge(c_1'/d)X_D/(D_md^{\,m})$. If
$c_4/D\ge2c_5\log(3D_m)/(D_md^{\,m})$, Proposition \ref{prop:energy} gives $\hh_\varphi(\alpha)\ge c_4/(2D)$.
Otherwise
$D_md^{\,m}<(2c_5/c_4)D\log3D$, and
$$\hh_\varphi(\alpha)\ \ge\ \frac{c_1'c_4}{2c_5d}\cdot\frac{X_D}{D\log3D}.$$
As $X_D\le1\le\log3D$ and $e_\p(P)\ge1$, both parts hold with the constant $c_2$ of Theorem \ref{thm:window} and
\[c_1=\tfrac12\min\bigl\{c_1'd^{-\nu-1},\ c_4/2,\ c_1'c_4/(2c_5d)\bigr\},\]
part (ii) with $D\log3D$ in place of $D$ and $c_2$ in place of $c_2+1$. Finally $\sqrt y\ge\log y$ for $y>0$, so
$\log3D\le\exp(\sqrt{\log3D})$, which gives (ii) as stated.
\end{proof}

\section{Points of bounded residue degree}\label{sec:resdeg}

In this section we prove  Theorem \ref{thm:resdeg}; $\Phi$ and $\varphi$ are the Latt\`es pair and map of \S\ref{sec:redp}. The proof
uses neither the growing congruence nor Theorem \ref{thm:window}. We use only that $[p]$ contracts the formal group $\p$-adically and
that the conjugates of $P$ lie in a bounded number of residue classes, and compare a single pair of levels of the orbit.

\subsection{The formal group}
We keep the embedding $\Kbar\hookrightarrow\CC_p$ and the valuation $v$ of \S\ref{sec:local}. We use the
following facts, valid over every closed subfield $\Omega\subseteq\CC_p$ containing $K_\p$, with $\mathfrak m_\Omega$ the maximal
ideal of its valuation ring. Let $\widehat E$ be the formal group of $E$, with parameter $t=-x/y$; its group law, its inverse
and the series $[m](T)$ have coefficients in $\OO$ \cite[IV.1--IV.2]{SilAEC}.
\begin{enumerate}[label=(F\arabic*),leftmargin=*]
\item Reduction is a surjective homomorphism $E(\CC_p)\to\widetilde E(\overline{\FF}_p)$. Its kernel $E^1$ consists of $O$ and the
points with $v(x)<0$, and $R\mapsto t(R)$ is an isomorphism of $E^1(\Omega)$ onto $\widehat E(\mathfrak m_\Omega)$; for $R\in E^1$,
$R\ne O$, one has $v(x(R))=-2v(t(R))$. This is \cite[VII.2.1, VII.2.2]{SilAEC}; surjectivity over $\CC_p$ follows from VII.2.1 applied over finite unramified extensions of $K_\p$, since every point of $\widetilde E(\overline{\FF}_p)$ is defined over a finite field.
\item For $R,R'\in E^1$, $v(t(R-R'))\ge\min\{v(t(R)),v(t(R'))\}$, since $t(R-R')$ is a power series in $t(R),t(R')$ with
coefficients in $\OO$ and without constant term.
\item For $m$ prime to $p$, $[m](T)=mT+\cdots$, so $[m]$ is an automorphism of $E^1(\Omega)$ \cite[IV.2.3]{SilAEC}.
\item Since $\widetilde E$ is supersingular, $\widehat E$ has height $2$ \cite[V.3.1]{SilAEC}, so $[p](T)\equiv h_0(T^{p^2})\pmod\p$ for
some $h_0\in k[[T]]$ \cite[IV.7.2]{SilAEC}. Thus $[p](T)=\sum_{i\ge1}a_iT^i$ with $a_1=p$, $a_i\in\p$ for $1\le i<p^2$, and $a_{p^2}\in\OO^\times$. Put
$$\varepsilon_E=\min_{1\le i<p^2}v(a_i),\qquad\frac1e\le\varepsilon_E\le1,\qquad \rho_E(s)=\min\{p^2s,\ s+\varepsilon_E\}.$$
Then $v\bigl(t([p]R)\bigr)\ge\rho_E\bigl(v(t(R))\bigr)$ for $R\in E^1$.
\end{enumerate}
For (F4): if $s=v(t(R))>0$, the terms with $i<p^2$ have valuation at least $\varepsilon_E+is\ge s+\varepsilon_E$, and the
others at least $p^2s$.

Each residue class of $E(\CC_p)$ contains a canonical torsion point of order prime to $p$, which plays the role of the periodic point of a residue disc.

\begin{lemma}\label{lem:teich}
Every residue class $R+E^1$ of $E(\CC_p)$ contains exactly one torsion point of order prime to $p$, written $\omega(\bar R)$. The map
$\bar R\mapsto\omega(\bar R)$ is a homomorphism, and if $R\in E(\Omega)$ then $\omega(\bar R)\in E(\Omega)$.
\end{lemma}

\begin{proof}
As $\widetilde E$ is supersingular, $\widetilde E(\overline{\FF}_p)$ has no point of order $p$ \cite[V.3.1]{SilAEC}, so $\bar R$ has order
$m$ prime to $p$. Then $[m]R\in E^1$, and by (F3) there is a unique $z\in E^1$ with $[m]z=[m]R$; the point $\omega=R-z$ satisfies
$[m]\omega=O$. Two torsion points of order prime to $p$ in the same class differ by a point of $E^1$ of order prime to $p$, which
is $O$ by (F3). If $R\in E(\Omega)$, then $[m]R\in E^1(\Omega)$ and $z\in E^1(\Omega)$ by (F3), so $\omega\in E(\Omega)$. The
homomorphism property follows from uniqueness.
\end{proof}

The next lemma quantifies the contraction of the formal group under $[p]$.

\begin{lemma}\label{lem:rho}
Let $e'\ge1$ and let $j_0\ge0$ be the least integer with $p^{2j_0}\ge\varepsilon_Ee'$. Put $u_0=1/e'$ and $u_{j+1}=\rho_E(u_j)$.
Then $(u_j)$ is increasing, and $u_j\ge\varepsilon_E(j-j_0+1)$ for $j\ge j_0$.
\end{lemma}

\begin{proof}
As $\rho_E(s)>s$ for $s>0$, $(u_j)$ is increasing. Let $a>0$ and $b=\min\{a,\varepsilon_E\}$. If $s\ge b$, then, as $\rho_E$ is
increasing,
$$\rho_E(s)\ \ge\ \rho_E(b)=\min\{p^2b,\ b+\varepsilon_E\}\ \ge\ \min\{p^2a,\ \varepsilon_E\};$$
by induction $u_j\ge\min\{p^{2j}/e',\varepsilon_E\}$, so $u_{j_0}\ge\varepsilon_E$. For $s\ge\varepsilon_E$ we have
$p^2s\ge s+\varepsilon_E$, so $\rho_E(s)=s+\varepsilon_E$, and the claim follows by induction.
\end{proof}

\subsection{Local heights}
For an embedding $\sigma$ of a number field $L'\supseteq K$ into some $\CC_\ell$ and $\beta,\beta'\in\mathbb P^1(L')$ put
$g_\sigma(\beta,\beta')=g_{\sigma|_K}(\sigma\beta,\sigma\beta')$, with $g_\iota$ as in \S\ref{sec:green}, and $\infty=(1:0)=x(O)$.

We collect the properties of the Green functions that are used in the proof.

\begin{lemma}\label{lem:pairlow}
\begin{enumerate}[label=\rm(\alph*)]
\item For $R\in E(L')$, $R\ne O$, $\frac12\sum_\sigma g_\sigma(x(R),\infty)=[L':\QQ]\,\hh_E(R)$.
\item For each embedding $\iota$ of $K$ into some $\CC_\ell$ there is $b_\iota\ge0$ with
$\bigl|\log\|\Phi^\iota(w)\|-d\log\|w\|\bigr|\le b_\iota$ for all $w\ne0$. Let $a_\iota=0$ if $\iota$ induces a place outside
$\mathcal S$, $a_\iota=\log2+2b_\iota/(d-1)+|r_\iota|$ if $\iota$ is archimedean, and $a_\iota=2b_\iota/(d-1)+|r_\iota|$ otherwise.
Then $g_\iota(z,z')\ge-a_\iota$ for all distinct $z,z'\in\mathbb P^1(\CC_\ell)$.
\item If $\sigma$ induces a place outside $\mathcal S$, then $g_\sigma(x(R),\infty)=\log^+|x(R)|_\sigma$, where
$\log^+=\max\{0,\log\}$. If $\sigma\colon L'\hookrightarrow\CC_p$ extends the fixed embedding of $K$ and $\sigma R\in E^1$, then
$\frac12g_\sigma(x(R),\infty)=v\bigl(t(\sigma R)\bigr)\log p$.
\end{enumerate}
Fixing the $b_\iota$, put $C_{\mathrm{bad}}=[K:\QQ]^{-1}\sum_\iota a_\iota$; it is finite, and $C_{\mathrm{bad}}\ge\log2$.
\end{lemma}

\begin{proof}
(a) Take the lifts $(x(R),1)$ and $(1,0)$, whose wedge is $-1$. Summing the definition over $\sigma$, the product formula and
$\Res(\Phi)\in K^\times$ give $\sum_\sigma r_\sigma=0$, and (G3) gives $[L':\QQ](\hh_\varphi(x(R))+\hh_\varphi(\infty))$. As
$\Phi(1,0)=(1,0)$, $\hh_\varphi(\infty)=0$, and $\hh_\varphi(x(R))=2\hh_E(R)$ by Lemma \ref{lem:hE}.
(b) An upper bound for $\log\|\Phi^\iota(w)\|-d\log\|w\|$ follows from the triangle inequality, and a lower bound from the identities expressing $\Res(\Phi)X^{2d-1}$ and $\Res(\Phi)Z^{2d-1}$ as combinations of $\Phi_0,\Phi_1$ with polynomial coefficients; this gives $b_\iota$. Outside $\mathcal S$ the claim is (G2). In general
$g_\iota=-\log\frac{|\tilde z\wedge\tilde z'|}{\|\tilde z\|\|\tilde z'\|}+(\widehat H_\iota-\log\|\cdot\|)(\tilde z)
+(\widehat H_\iota-\log\|\cdot\|)(\tilde z')-r_\iota$. The first term is $\ge-\log2$ if $\iota$ is archimedean and $\ge0$
otherwise, and with $G=\Phi^\iota$,
$$\widehat H_\iota(w)-\log\|w\|=\sum_{j\ge0}d^{-j-1}\Bigl(\log\|G^{(j+1)}(w)\|-d\log\|G^{(j)}(w)\|\Bigr)$$
has absolute value at most $b_\iota/(d-1)$. Finally $C_{\mathrm{bad}}$ is finite as $\mathcal S$ is, and $C_{\mathrm{bad}}\ge\log2$
as there are $[K:\QQ]$ archimedean embeddings.
(c) By (G2), $g_\sigma(x(R),\infty)=-\log\bigl(1/\max\{|x(R)|_\sigma,1\}\bigr)=\log^+|x(R)|_\sigma$. If $\sigma$ extends the fixed
embedding, then $|x(R)|_\sigma=p^{-v(x(\sigma R))}$, and (F1) gives $\log^+|x(R)|_\sigma=2v(t(\sigma R))\log p$ for $\sigma R\in E^1$.
\end{proof}

\subsection{Proof of Theorem \ref{thm:resdeg}}
Let $P$ be non-torsion, $D'=[K(P):K]$, $P_1,\dots,P_{D'}$ the $K$-conjugates of $P$, viewed in $E(\CC_p)$, $e'=e_\p(P)$, and
$L'\subset\Kbar$ a finite Galois extension of $K$ over which all $P_i$ are defined.

Every prime of $K(P)$ above $\p$ has residue degree at least $f$, so we may assume $f_0\ge f$. Let $U\subset\widetilde E(\overline{\FF}_p)$ be the union of the groups $\widetilde E(\FF_{p^j})$ over the $j\le f_0$
divisible by $f$, let $N_0=|U|$, let $\mu_0$ be the exponent of the finite group generated by $U$, which is prime to $p$, and let
$\kappa_0\ge1$ be the order of $p$ modulo $\mu_0$. All these depend only on $E$, $\p$ and $f_0$. Let $\Omega_i$ be the closure of
$K_\p(P_i)$ in $\CC_p$; it is the completion of $K(P_i)\cong K(P)$ at a prime above $\p$, so its absolute ramification index
$e_i$ is at most $e'$, and its residue field contains $k$ and has degree at most $f_0$ over $\FF_p$. Hence $\bar P_i\in U$.

Put $\omega_i=\omega(\bar P_i)$ and $z_i=P_i-\omega_i$. By Lemma \ref{lem:teich}, $\omega_i\in E(\Omega_i)$, so
$z_i\in E^1(\Omega_i)$, and $z_i\ne O$ as $P$ is not torsion. Thus $t(z_i)$ is a non-zero element of $\Omega_i$ of positive
valuation, so $v(t(z_i))\ge1/e_i\ge1/e'$. By (F4) and since $\rho_E$ is increasing, $v\bigl(t([p^j]z_i)\bigr)\ge u_j$ for all
$j\ge0$, with $u_j$ as in Lemma \ref{lem:rho}.

Fix $n\ge0$, put $m=n+\kappa_0$ and $R_{ii'}=[p^n]P_i-[p^m]P_{i'}\in E(L')$. If $R_{ii'}=O$, then $[p^n]P$ and
$[p^m]P$ are conjugate over $K$, so $p^{2n}\hh_E(P)=p^{2m}\hh_E(P)$ and $P$ would be torsion; hence $R_{ii'}\ne O$. By
Lemma \ref{lem:pairlow}(a) and the parallelogram law,
\begin{equation}\label{eq:resdegE}
\mathcal E=\frac12\sum_\sigma\sum_{i,i'}g_\sigma\bigl(x(R_{ii'}),\infty\bigr)=[L':\QQ]\sum_{i,i'}\hh_E(R_{ii'})
\ \le\ 2[L':\QQ]\,D'^2\,(p^{2n}+p^{2m})\,\hh_E(P).
\end{equation}

By Lemma \ref{lem:pairlow}(b), the $\sigma$ inducing places of $\mathcal S$ contribute at
least $-\frac12D'^2[L':K]\sum_\iota a_\iota=-\frac12[L':\QQ]D'^2C_{\mathrm{bad}}$, and by (G2) the other $\sigma$ contribute
non-negative terms. Let $\sigma\colon L'\hookrightarrow\CC_p$ induce $\p$ on $K$; there are $[L':K]\,[K_\p:\QQ_p]=[L':K]\,ef$ of them.
The restriction of $\sigma$ to $K$ and the given embedding of $K$ both extend to $\QQ_p$-embeddings of $K_\p$, which differ
by an element of $\mathrm{Gal}(\overline\QQ_p/\QQ_p)$; extending it by continuity gives an automorphism $\theta$ of $\CC_p$ with
$v\circ\theta=v$ and $\sigma=\theta\circ\gamma$ on $L'$ for some $\gamma\in\mathrm{Gal}(L'/K)$. Then $|x(R)|_\sigma=|\theta(x(\gamma R))|_p=|x(\gamma R)|_p$, and $\gamma R_{ii'}=R_{\gamma i,\gamma i'}$, where $\gamma P_i=P_{\gamma i}$.
Hence, by Lemma \ref{lem:pairlow}(c),
$$\frac12\sum_{i,i'}g_\sigma\bigl(x(R_{ii'}),\infty\bigr)=\frac12\sum_{i,i'}\log^+\bigl|x(R_{ii'})\bigr|_p,$$
where on the right $R_{ii'}$ is viewed in $E(\CC_p)$ through the fixed embedding.
Suppose $\bar P_i=\bar P_{i'}$, so that $\omega_i=\omega_{i'}=:\omega$. As $\omega$ lies in the group $\omega(G)$, where $G$ is generated by $U$, $[p^{\kappa_0}]\omega=\omega$, and
therefore $R_{ii'}=[p^n]z_i-[p^m]z_{i'}\in E^1$. By (F2) and Lemma \ref{lem:rho},
$v(t(R_{ii'}))\ge\min\{u_n,u_m\}=u_n$, so the term of $(i,i')$ is at least $u_n\log p$ by Lemma \ref{lem:pairlow}(c). The
other terms are non-negative. Since the $\bar P_i$ lie in $U$, the number of ordered pairs with $\bar P_i=\bar P_{i'}$ is at
least $D'^2/N_0$ by the Cauchy--Schwarz inequality. Hence
$$\mathcal E\ \ge\ [L':K]\,ef\log p\cdot\frac{D'^2}{N_0}\,u_n\ -\ \frac12[L':\QQ]\,D'^2\,C_{\mathrm{bad}}.$$

With $\Lambda=N_0[K:\QQ]/(ef\log p)$, comparison with \eqref{eq:resdegE} gives, after division by $[L':\QQ]D'^2$,
$$\frac{u_n}{\Lambda}\ \le\ \frac{C_{\mathrm{bad}}}2+4p^{2(n+\kappa_0)}\,\hh_E(P)\qquad\text{for every }n\ge0.$$
Take $n=j_0-1+\lceil \Lambda C_{\mathrm{bad}}/\varepsilon_E\rceil$, with $j_0$ as in Lemma \ref{lem:rho}; then $n\ge j_0$, and
$u_n\ge\varepsilon_E\lceil \Lambda C_{\mathrm{bad}}/\varepsilon_E\rceil\ge \Lambda C_{\mathrm{bad}}$, so $4p^{2(n+\kappa_0)}\hh_E(P)\ge C_{\mathrm{bad}}/2$.
By minimality of $j_0$, $p^{2(j_0-1)}<\varepsilon_Ee'\le e'$ if $j_0\ge1$, and $p^{2(j_0-1)}<1\le e'$ if $j_0=0$. Using
$\varepsilon_E\ge1/e$,
$$\hh_E(P)\ \ge\ \frac{C_{\mathrm{bad}}}{8\,p^{2\kappa_0+2\lceil e\Lambda C_{\mathrm{bad}}\rceil}}\cdot\frac1{e'}=:\frac{c_3}{e_\p(P)}.$$
Finally $e_\p(P)\le e[K(P):K]\le2eD$. \qed

\begin{remark}[Sharpness]\label{rem:resdegsharp}
The exponent of $e_\p(P)$ in Theorem \ref{thm:resdeg} cannot be lowered. Suppose $Q_0\in E(K)\cap E^1$ is non-torsion with
$v_\p(t(Q_0))=1$, and let $n\ge1$. By (F4), $[p](T)\equiv h_0(T^{p^2})\pmod\p$ with $h_0(T)=\bar a_{p^2}T+\cdots$, $\bar a_{p^2}\ne0$; iterating,
$[p^n](T)$ has coefficients in $\p$ below degree $p^{2n}$ and a unit coefficient in degree $p^{2n}$, so by the Weierstrass preparation theorem $[p^n](T)-t(Q_0)$ is a unit times an Eisenstein polynomial of degree
$p^{2n}$ over $\OO_{K_\p}$. Its roots are the values $t(P)$ for the $p^{2n}$ points of $[p^n]^{-1}(Q_0)$, all of which lie in
$E^1$. Hence every $P\in[p^n]^{-1}(Q_0)$ has $[K(P):K]=p^{2n}$, a single prime above $\p$, totally ramified over $\p$, so that
$e_\p(P)=ep^{2n}$ and the residue degree is $f$, and $\hh_E(P)=\hh_E(Q_0)p^{-2n}=e\,\hh_E(Q_0)/e_\p(P)$. Here also $p^{2n}/2\le D\le p^{2n}$, so $\hh_E(P)\le\hh_E(Q_0)/D$ and the exponent $1$ in Theorem \ref{thm:main} is optimal as well.
Such $Q_0$ exist: for $E\colon y^2+y=x^3-x$ over $\QQ$, which is supersingular for instance at $p=3$, $17$ and $19$, with
$n_p=|\widetilde E(\FF_p)|=7$, $18$, $20$, the point $Q_0=[n_p](0,0)$ is non-torsion and satisfies $v_p(t(Q_0))=1$ in each case.
\end{remark}

\section{Supersingular primes over number fields}\label{sec:ss}

In this section we examples for which a prime of good supersingular reduction is known to exist; Theorems \ref{thm:main}--\ref{thm:abelian} only need one of those. 

Such primes exist for curves with complex multiplication and for curves whose $j$-invariant generates a field with a real embedding.

\begin{proposition}\label{prop:ss}
Let $E$ be an elliptic curve over a number field $K$. Then $E$ has infinitely many primes of good supersingular reduction in
each of the following cases:
\begin{enumerate}[label=\rm(\alph*)]
\item $E$ has complex multiplication;
\item the field $\QQ(j(E))$ has a real embedding; in particular, $E$ is defined over a field with a real embedding, or
$j(E)\in\QQ$.
\end{enumerate}
\end{proposition}

\begin{proof}
Whether the reduction at a prime $\mathfrak Q$ of good reduction is supersingular depends only on $j(E)\bmod\mathfrak Q$, since
supersingularity is a property of the reduced curve over an algebraic closure of the residue field \cite[V.3.1]{SilAEC}. Also,
$E$ has bad reduction at only finitely many primes.

(a) Let $E$ have complex multiplication by an order in an imaginary quadratic field $F_0$. By Deuring's criterion
\cite[Ch.~13]{Lan}, at a prime of good reduction above a rational prime that is inert in $F_0$ the reduction is
supersingular. There are infinitely many such rational primes.

(b) Let $K_0=\QQ(j(E))$. There is an elliptic curve $E_0/K_0$ with $j(E_0)=j(E)$ \cite[III.1.4]{SilAEC}. Since $K_0$ has a
real embedding, Elkies \cite{Elk89} shows that $E_0$ has infinitely many primes of supersingular reduction. At every prime
$\mathfrak Q$ of $K$ above such a prime at which $E$ and $E_0$ both have good reduction, $j(E)\bmod\mathfrak Q$ is a supersingular
$j$-invariant, so $E$ has supersingular reduction at $\mathfrak Q$.
\end{proof}

For $E$ without complex multiplication and $\QQ(j(E))$ totally complex, the existence of a supersingular prime is not known in
general; in specific cases one can often be found by computation. For example, over $K=\QQ(i)$ the curves $y^2+y=x^3+ix$ and $y^2=x^3+ix+1$ have
$j$-invariants $110592(64-27i)/4825$ and $-6912i(27+4i)/745$, which are neither real nor integral, and they have good
supersingular reduction at $(1+i)$ and at $(3)$ respectively.

\section{Abelian varieties isogenous to products}\label{sec:abelian}

In this section we prove Theorem \ref{thm:abelian} and give examples of abelian varieties to which it applies.

Let $A$ be an abelian variety over $K$. We use the following properties of the N\'eron--Tate height $\hh_{A,\mathcal L}$ of an ample
symmetric line bundle $\mathcal L$ \cite[\S B.5]{HSiBook}: it is non-negative and vanishes exactly at the torsion points; it is
additive in $\mathcal L$ on symmetric line bundles; $\hh_{B,\mathcal M}(\phi(P))=\hh_{A,\phi^*\mathcal M}(P)$ for a homomorphism
$\phi\colon A\to B$; and for $B=E_1\times\dots\times E_g$ and $\mathcal M=\bigotimes_i\mathrm{pr}_i^*\OO_{E_i}(O_i)$,
$$\hh_{B,\mathcal M}(Q_1,\dots,Q_g)\ =\ \sum_{i=1}^g\hh_{E_i}(Q_i),$$
where $\hh_{E_i}$ is the N\'eron--Tate height of $\OO_{E_i}(O_i)$, which is the height used in the rest of the paper. Since these
heights on $A(\Kbar)$ do not depend on the field of definition, we may compute them after extending the base field to $K'$.

The next lemma compares the N\'eron--Tate heights of two ample symmetric line bundles.

\begin{lemma}\label{lem:ample}
Let $\mathcal L,\mathcal L'$ be ample symmetric line bundles on $A$. There is $n\ge1$ with $\hh_{A,\mathcal L'}\le n\,\hh_{A,\mathcal L}$ on
$A(\Kbar)$.
\end{lemma}

\begin{proof}
For $n$ large, $\mathcal L^{\otimes n}\otimes\mathcal L'^{-1}$ is ample, and it is symmetric. Its N\'eron--Tate height is therefore
non-negative, and by additivity it equals $n\hh_{A,\mathcal L}-\hh_{A,\mathcal L'}$.
\end{proof}

\begin{proof}[Proof of Theorem \ref{thm:abelian}]
Replacing $K'$ by a Galois closure of $K'/K$, we may assume that $K'/K$ is Galois: the isogeny stays defined, and good
supersingular reduction is preserved at the primes above the chosen ones. Let $\phi\colon A_{K'}\to B=E_1\times\dots\times E_g$ be the
isogeny and $\mathcal M$ as above. Then $\phi^*\mathcal M$ is ample, as $\phi$ is finite, and symmetric, as $\phi$ commutes with $[-1]$. By
Lemma \ref{lem:ample} there is $n$ with $\hh_{A,\phi^*\mathcal M}\le n\,\hh_{A,\mathcal L}$, so
$$\hh_{A,\mathcal L}(P)\ \ge\ \frac1n\hh_{B,\mathcal M}(\phi(P))\ =\ \frac1n\sum_{i=1}^g\hh_{E_i}(Q_i),\qquad(Q_1,\dots,Q_g)=\phi(P).$$
If $P$ is not torsion, neither is $\phi(P)$, since $\ker\phi$ is finite; so some $Q_i$ is not torsion. As $\phi$ is defined over
$K'$, the point $Q_i$ is defined over $K'(P)$, so $D_i=[K'(x(Q_i)):K']$ divides $[K'(P):K']$. Since $K'/K$ is Galois,
$[K'(P):K']=[K(P):K(P)\cap K']$, which divides $D$. Hence $D_i\le D$. Theorem
\ref{thm:main}(ii) for $E_i$ over $K'$ gives $\hh_{E_i}(Q_i)\ge c_{1,i}D_i^{-1}\exp(-(c_{2,i}+1)\sqrt{\log3D_i})$, and the right
side decreases with $D_i$; so the claim holds with $c=\min_ic_{1,i}/n$ and $c'=\max_ic_{2,i}+1$.
\end{proof}

The hypothesis of Theorem \ref{thm:abelian} holds in the following cases.

\begin{corollary}\label{cor:abex}
The conclusion of Theorem \ref{thm:abelian} holds in each of the following cases:
\begin{enumerate}[label=\rm(\alph*)]
\item $A$ is isogenous over $K$ to a product of elliptic curves $E_i$ over $K$, each with complex multiplication or with
$\QQ(j(E_i))$ having a real embedding; in particular, every abelian variety over $\QQ$ that is isogenous over $\QQ$ to a
product of elliptic curves;
\item $A=\mathrm{Res}_{F/K}E$ is the Weil restriction of an elliptic curve $E$ over a finite extension $F/K$, where $E$ has
complex multiplication or $\QQ(j(E))$ has a real embedding.
\end{enumerate}
If $K=\QQ$, $F$ is quadratic and $E$ is not isogenous over $F$ to the base change of an elliptic curve over $\QQ$, then
$\mathrm{Res}_{F/\QQ}E$ in (b) is a simple abelian surface over $\QQ$.
\end{corollary}

\begin{proof}
(a) Each $E_i$ has a supersingular prime by Proposition \ref{prop:ss}; take $K'=K$.

(b) Let $K'$ be a Galois closure of $F/K$. Over $K'$, $\mathrm{Res}_{F/K}E$ is isomorphic to $\prod_\sigma E^\sigma$, the product over the
$K$-embeddings $\sigma\colon F\to K'$ \cite[\S1]{Mil}. Each $E^\sigma$ has complex multiplication if $E$ does, and
$\QQ(j(E^\sigma))=\sigma(\QQ(j(E)))$ is isomorphic to $\QQ(j(E))$, so it has a real embedding if $\QQ(j(E))$ does. By Proposition
\ref{prop:ss} each $E^\sigma$ has a supersingular prime, and Theorem \ref{thm:abelian} applies.

For the last assertion, $A=\mathrm{Res}_{F/\QQ}E$ is an abelian surface. If it were not simple over $\QQ$, it would be isogenous over
$\QQ$ to $E_1\times E_2$ with elliptic curves $E_i$ over $\QQ$. Over $F$, $A_F\cong E\times E^\sigma$, where $\sigma$ is the non-trivial
automorphism of $F$, so $E\times E^\sigma$ would be isogenous over $F$ to $(E_1)_F\times(E_2)_F$, and by Poincar\'e's complete
reducibility theorem $E$ would be isogenous over $F$ to $(E_1)_F$ or $(E_2)_F$.
\end{proof}

\end{document}